\documentclass[11pt,letterpaper]{article}
\usepackage[letterpaper,margin=1.05in]{geometry}
\usepackage{amsmath,amssymb,amsthm}
\usepackage{microtype}
\usepackage{tikz}
\usetikzlibrary{decorations.pathreplacing}
\usepackage[title,titletoc]{appendix}
\usepackage{etoolbox}
\usepackage[bottom]{footmisc}
\usepackage[hidelinks,bookmarksnumbered,hyperfootnotes=false]{hyperref}
\hypersetup{pdftitle={Mixing time and logarithmic Sobolev inequalities for the torus shuffle},
  pdfauthor={Chen Shang},
  pdfsubject={Probability, mixing times, logarithmic Sobolev inequalities}}
\newcommand{\paperdate}{November 22, 2026}
\newcommand{\paperemail}{cs3642@columbia.edu}
\makeatletter
\patchcmd{\l@section}{\addvspace{1.0em \@plus\p@}}{\addvspace{0.5em \@plus\p@}}{}{}
\patchcmd{\section}{\normalfont\Large\bfseries}{\normalfont\Large\bfseries\raggedright}{}{}
\patchcmd{\thebibliography}{\section*{\refname}}
  {\section*{\refname}\phantomsection\addcontentsline{toc}{section}{\refname}}
  {}{\errmessage{The patch of thebibliography failed}}
\makeatother
\newtheorem{theorem}{Theorem}[section]
\newtheorem{lemma}[theorem]{Lemma}
\newtheorem{proposition}[theorem]{Proposition}

\theoremstyle{remark}
\newtheorem*{remark}{Remark}
\numberwithin{equation}{section}
\newcommand{\Ent}{\operatorname{Ent}}
\newcommand{\Var}{\operatorname{Var}}
\newcommand{\E}{\mathbb E}
\newcommand{\one}{\mathbf{1}}
\newcommand{\cE}{\mathcal E}
\newcommand{\tmix}{t_{\mathrm{mix}}}
\newcommand{\step}[1]{\par\smallskip\noindent\emph{#1}\enspace}
\allowdisplaybreaks[1]
\begin{document}

\begin{center}
{\LARGE\bfseries Mixing time and logarithmic Sobolev inequalities\\[3pt]
for the torus shuffle\par}
\vspace{14pt}
{\large Chen Shang\par}
\end{center}
\begingroup
\renewcommand{\thefootnote}{}
\footnotetext{\emph{Date}: \paperdate.}
\footnotetext{Email: \href{mailto:\paperemail}{\paperemail}.}
\endgroup
\begin{abstract}
The torus shuffle acts on an $n\times n$ board of $n^2$ cards: each
step chooses a row or a column uniformly at random and rotates it
cyclically by one unit in a uniformly chosen direction. We prove that
the half-lazy version of this chain has total variation mixing time
of order $n^3\log n$ and logarithmic Sobolev constant of order
$n^{-3}$ on each communicating class. The order of the mixing time is
the one conjectured by Diaconis. The proof of the upper bound has two
steps. First, the entropy of a function on the board is at most a
constant times $n$ times its mean loss when a single row or column is
completely rearranged. Second, this mean loss is at most a constant
times $n^2$ times the Dirichlet form of the shuffle. The second step
bisects each line recursively and charges the cost of a column
rotation only to the intervals that contain the column, so that the
costs over all scales add up without a logarithmic factor. A Fourier
observable gives the matching lower bound on the mixing time.

\medskip
\noindent\emph{2020 Mathematics Subject Classification}: Primary 60J10;
Secondary 60B15, 39B62.

\noindent\emph{Keywords}: card shuffling, entropy factorization,
logarithmic Sobolev inequality, mixing time, torus shuffle.
\end{abstract}

\setcounter{tocdepth}{1}
\tableofcontents

\section{Introduction}\label{sec:intro}

Place $n^2$ distinct cards on the sites of an $n\times n$ board with
periodic boundary conditions. A \emph{line} is a row or a column of
the board. One move of the \emph{torus shuffle} chooses one of the
$2n$ lines uniformly at random and rotates it cyclically by one site,
in a direction chosen uniformly at random (Figure~\ref{fig:move}). We
study the half-lazy chain: at each update it stays fixed with
probability $1/2$ and otherwise performs such a move. Laziness makes
the spectrum of the chain nonnegative, which Section~\ref{sec:mixing}
uses to pass from the logarithmic Sobolev inequality to the decay of
entropy in discrete time. It also removes a periodicity: a unit
rotation is a cycle of length $n$, so for even $n$ every move is an
odd permutation of the cards, and the shuffle without holding has
period two. Blumberg, Morris and Senda use the same holding
probability \cite{bms}.

\begin{figure}[ht]
\centering
\begin{tikzpicture}[x=0.6cm,y=0.6cm,>=stealth]
  \fill[black!14] (0,2) rectangle (5,3);
  \draw[black!45,step=1] (0,0) grid (5,5);
  \foreach \r in {0,1,2,3,4}
    \foreach \c in {0,1,2,3,4} {
      \pgfmathtruncatemacro{\k}{5*\r+\c+1}
      \node at (\c+0.5,4.5-\r) {\small\k};
    }
  \draw[->,thick] (5.8,2.5) -- (8.2,2.5);
  \begin{scope}[shift={(9,0)}]
    \fill[black!14] (0,2) rectangle (5,3);
    \draw[black!45,step=1] (0,0) grid (5,5);
    \foreach \r in {0,1,3,4}
      \foreach \c in {0,1,2,3,4} {
        \pgfmathtruncatemacro{\k}{5*\r+\c+1}
        \node at (\c+0.5,4.5-\r) {\small\k};
      }
    \foreach \c/\k in {0/15,1/11,2/12,3/13,4/14}
      \node at (\c+0.5,2.5) {\small\k};
  \end{scope}
\end{tikzpicture}
\caption{One move of the torus shuffle on a $5\times5$ board. The
shaded row is rotated by one site to the right. The card at its right
end, here card~$15$, re-enters at the left end, and the cards outside
the row do not move.}
\label{fig:move}
\end{figure}

Diaconis listed the shuffle without holding among open problems on
random walks and
conjectured that order $n^3\log n$ moves suffice to mix the board
\cite[Ch.~4, \S G, Problem~4]{diaconis}. For the half-lazy chain,
Blumberg, Morris and Senda obtained an $O(n^3\log^3 n)$ upper bound
on the mixing time through an entropy argument based on collisions of
triples of cards \cite{bms}. We prove that the mixing time of this
chain has the conjectured order $n^3\log n$, and that the chain
satisfies a logarithmic Sobolev inequality with constant of order
$n^{-3}$. The second remark after Theorem~\ref{thm:main} says what
this gives for the shuffle without holding.

The order $n^3$ is the time scale of a single card, as Diaconis
observed \cite{diaconis}. A card changes
its horizontal coordinate with probability of order $1/n$ per update,
so its horizontal motion on a cycle of length $n$ has relaxation time
of order $n^3$. The motions of different cards are dependent,
however, since each move translates an entire line. To bound the
mixing time of the whole board we argue in two steps. We first bound
the entropy of a function on the board by its losses under complete
rearrangements of single rows and columns. We then bound these losses
by the Dirichlet form of the shuffle, building local rearrangements
out of unit rotations.
Subsection~\ref{subsec:outline} states the two estimates and explains
how they give the upper bound. It ends with a sketch of the proof of
the matching bounds.

\subsection{Model and main result}\label{subsec:model}
Let $V=(\mathbb Z/n\mathbb Z)^2$ be the set of sites and put $N=n^2$.
A \emph{board} is a bijection $x$ from a set of $N$ distinct labels
to $V$; it records the site of each card. The transition kernel $P$
of the half-lazy torus shuffle has $P(x,x)=1/2$, and each of the $4n$
signed unit rotations, that is, a line together with a direction, is
applied with probability $1/(8n)$. When two signed rotations coincide
as permutations, as happens for $n=2$, where the two directions
agree, their probabilities add.

The rotations permute the sites, and the cards move with them. Let
$G_n$ be the group of permutations of $V$ generated by the unit
rotations. The chain started at a board $x$ stays in the orbit $G_nx$,
and its stationary law there is $\pi_x=\mathrm{Unif}(G_nx)$.
Proposition~\ref{prop:classes} shows that $G_n=S_N$ for even
$n$ and $G_n=A_N$ for odd $n$. Thus the chain has one communicating
class when $n$ is even and two classes, distinguished by parity, when
$n$ is odd.

Fix a communicating class and write $\pi$ for its uniform law. Unless
otherwise specified, all norms and inner products are taken with
respect to this probability measure. We use the total variation
convention $\|\mu-\nu\|_{\rm TV}=\frac12\sum_x|\mu(x)-\nu(x)|$.
For $h\geq0$ and real $f$ put
\[
 \Ent_\pi(h)=\E_\pi[h\log h]-\E_\pi[h]\log\E_\pi[h],\qquad
 \cE_P(f,f)=\langle f,(I-P)f\rangle_\pi ,
\]
where $\E_\pi$ is the expectation under $\pi$, with $0\log0=0$ and
natural logarithms. The \emph{logarithmic Sobolev
constant} and the \emph{mixing time} of the chain are
\[
 \alpha(P)=\inf_{\Ent_\pi(f^2)>0}
       \frac{\cE_P(f,f)}{\Ent_\pi(f^2)},\qquad
 \tmix^P(\delta)=\min\left\{t\geq0:
    \sup_x\|P^t(x,\cdot)-\pi_x\|_{\rm TV}\leq\delta\right\}.
\]
The infimum runs over functions on the fixed class. Relabeling the
cards maps one class onto the other and commutes with the dynamics,
so $\alpha(P)$ does not depend on the class. See \cite{dsc96} for
logarithmic Sobolev inequalities on finite Markov chains.

\begin{theorem}\label{thm:main}
There are absolute constants $c,C>0$ such that for every $n\geq2$,
\[
 c n^{-3}\leq\alpha(P)\leq C n^{-3},
 \qquad c n^3\log n\leq\tmix^P(1/4)\leq C n^3\log n .
\]
For $0<\delta\leq1/4$ one has
$\tmix^P(\delta)\leq C n^3(\log n+\log(1/\delta))$.
\end{theorem}

\begin{remark}
The lower bound on $\alpha(P)$ and the upper bounds on the mixing
time come from the two comparison estimates described in
Subsection~\ref{subsec:outline}. The upper bound on $\alpha(P)$ and
the lower bound on the mixing time come from one Fourier
eigenfunction, the latter by the eigenfunction method of
Wilson~\cite{wilson}. The constants $c$ and $C$ are not explicit:
several estimates in Section~\ref{glob:section} are proved for all
sufficiently large $n$, and the remaining board sizes are handled by
a compactness argument.
\end{remark}

\begin{remark}[The shuffle without holding]
Let $\widetilde P$ be the kernel of the shuffle without the holding
step, which is the chain considered in
\cite[Ch.~4, \S G, Problem~4]{diaconis}. Then $P=(I+\widetilde P)/2$
and $\cE_{\widetilde P}(f,f)=2\cE_P(f,f)$ for every $f$. Hence
$\alpha(\widetilde P)=2\alpha(P)$, and the bounds of
Theorem~\ref{thm:main} on the logarithmic Sobolev constant hold for
$\widetilde P$ with the constants doubled. For even $n$ the kernel
$\widetilde P$ has period two, so the distance to equilibrium of its
chain does not tend to zero. For odd $n$ the kernel $\widetilde P$ is
aperiodic. Theorem~\ref{thm:main} does not bound the mixing time of
$\widetilde P$ in discrete time, because Step~2 of
Subsection~\ref{subsec:decay} uses that the spectrum of $P$ is
nonnegative.
\end{remark}

\subsection{Related work}\label{subsec:related}

\paragraph{The model.}
Diaconis introduced the shuffle as a variant of the fifteen puzzle
\cite{diaconis}. On the $n\times n$ torus that puzzle has mixing time
of order $n^4\log n$ \cite{morrisraymer}. Blumberg, Morris and Senda
extend to
triples of cards the entropy method that Morris introduced for pairs
\cite{bms,morris}. They point out that comparison with random
transpositions \cite{dsc93} gives only $O(n^4\log n)$.

Each move of the torus shuffle is a cycle of length $n$ on $N=n^2$
cards. For the walk that applies at each step a uniformly chosen
cycle of length $k=o(N)$, the mixing time is asymptotic to
$(N/k)\log N$. Berestycki, Schramm and Zeitouni proved this for fixed
$k$ \cite{bsz}, and Hough proved it in the whole range \cite{hough}.
Nestoridi and Olesker-Taylor found the limit profile of the distance
to equilibrium in the same range \cite{nestoridiolesker}. For
$k=o(N/(\log N)^4)$ the
distribution at the first time when every card has been moved is
described in \cite{sszz}. For $k=n$ the mixing time is of order
$n\log n$. Restricting the cycles to unit rotations of lines slows
mixing by a factor of order $n^2$, the diffusive time scale of a
nearest-neighbor walk on a line of $n$ sites.

\paragraph{Entropy factorization.}
Our first step bounds the entropy of a function by its losses when
the cards of single lines are rearranged uniformly. Bristiel and
Caputo study inequalities of this form for block shuffles, in which
the cards on a random block of sites are rearranged uniformly
\cite{bristielcaputo}. They find the optimal constants when the rate
of a block depends only on its size, and they ask when the constant
for the whole deck equals the constant for functions of a single
card. For even $n$ the rearrangements of lines are block shuffles.
For arbitrary weighted blocks, Caputo and Salez bound the constant
from below by the least total weight of the blocks that contain a
given pair of sites \cite[Theorem~4]{caputosalez}. No line contains
two sites in different rows and different columns, so this bound is
zero for the lines. Theorem~\ref{thm:global} shows that for the lines
the constant for the whole deck is within an absolute factor of the
constant for a single card. Its proof rests on the bilinear
inequality \eqref{glob:bilinear}, an analog for rows and columns of
the inequality of Carlen, Lieb and Loss for functions of single cards
\cite{cll}.

Factorizations of entropy over blocks are also a standard
tool for spin systems \cite{caputoparisi}. A block shuffle on a line,
in which a random block of $k$ consecutive cards is rearranged, was
introduced by Nestoridi, Priestley and Schmid \cite{nestoridips}. They
find the order of its mixing time when $k=o(m^{2/3})$, for a deck of
$m$ cards. Rearrangements of short intervals of a line enter the
proof of Theorem~\ref{thm:physical}.

\paragraph{Single cards and the lower bound.}
For interchange processes, in which a move transposes the two cards
on an edge of a graph, the spectral gap of the deck equals that of a
single card. Aldous conjectured this, and Caputo, Liggett and
Richthammer proved it \cite{clr}. For $n\geq3$ a unit rotation moves
more than two cards, so the torus shuffle is not an interchange
process. Theorem~\ref{thm:main} shows that its logarithmic Sobolev
constant is of the order of the inverse time scale of a single card.
For the lower bound on the mixing time we use the eigenfunction
method of Wilson \cite{wilson}, which is also presented in
\cite[Theorem~13.28]{lpw}. For the adjacent transposition shuffle of
$n$ cards the method gives a lower bound of order $n^3\log n$. Lacoin
proved that this bound is sharp, so the shuffle has a cutoff
\cite{lacoin}. Nam and Nestoridi proved a cutoff for a systematic scan
version of this shuffle, with a form of the method that allows
approximate eigenfunctions \cite{namnestoridi}.

\subsection{Outline of the proof}\label{subsec:outline}
We compare the shuffle with the dynamics that rearranges all the
cards on a line at once. For a line $\ell$, the \emph{line group} of
$\ell$ consists of the permutations of the sites of $\ell$ that
preserve the communicating class, all other sites being fixed. It is
the symmetric group on the $n$ sites of $\ell$ for even $n$ and the
alternating group for odd $n$. Let $B_\ell$ denote the average over
the line group: $B_\ell h(x)$ is the mean of $h(\sigma x)$ over
$\sigma$ in the line group of $\ell$. This is the heat bath that
resamples the arrangement of the cards on $\ell$. Define the mean
entropy loss
\[
 D_{\rm line}(h)=\frac{1}{2n}\sum_\ell
 \bigl[\Ent_\pi(h)-\Ent_\pi(B_\ell h)\bigr],
\]
the expected loss of entropy when the cards on a uniformly chosen
line are rearranged. The upper bounds on the mixing time in
Theorem~\ref{thm:main} rest
on two estimates, valid for every $h\ge0$ with absolute constants
$C_1,C_2$:
\begin{equation}\label{eq:bridge}
 \Ent_\pi(h)\leq C_1 nD_{\rm line}(h),\qquad
 D_{\rm line}(h)\leq C_2 n^2\cE_P(\sqrt h,\sqrt h).
\end{equation}
The first is part of Theorem~\ref{thm:global} and the second is
Theorem~\ref{thm:physical}. Combining the estimates with $h=f^2$ gives
$\Ent_\pi(f^2)\leq Cn^3\cE_P(f,f)$, and hence the lower bound on
$\alpha(P)$. For the half-lazy chain this inequality makes entropy
contract by a factor $1-cn^{-3}$ per update. A point mass has entropy
$\log|G_n|\leq 2n^2\log n$, so $O(n^3\log n)$ updates reduce it to
a fixed constant. Pinsker's inequality then gives the upper bound on
the mixing time, and Subsection~\ref{subsec:decay} contains the
argument in discrete time. The bounds in the other direction are sketched at the end
of this outline and proved in Subsection~\ref{subsec:lower}.

\paragraph{The two estimates for a single card.}
Both estimates have a simple meaning when $h$ depends only on the site
of one card. The row and the column of a card are independent under
$\pi$, so the tensorization of entropy gives
$\Ent_\pi(h)\leq2nD_{\rm line}(h)$. The factor $n$ appears because a
uniformly chosen line contains the card with probability $1/n$. The
second estimate becomes a logarithmic Sobolev inequality for the
nearest-neighbor walk on a cycle of length $n$, whose diffusive time
scale is $n^2$. The product $n\cdot n^2$ is the time scale of a single
card from the introduction. Theorems~\ref{thm:global}
and~\ref{thm:physical} extend the two estimates, with absolute
constants, to functions of all the cards. There the rows and columns
of different cards are dependent, and a unit rotation moves a whole
line.

\paragraph{Factorization across rows and columns.}
The first estimate in \eqref{eq:bridge} concerns the two operations
that rearrange every row at once and every column at once.
Rearranging every row preserves the set of labels in each row, which
we call the row data of the board; rearranging every column preserves
the column data. The row data and the column data are dependent,
since every site must contain exactly one card. We prove an $L^p$
correlation inequality between them for some $p<2$ that does not
depend on $n$. For $p=2$ it is the Cauchy--Schwarz inequality, and
for $p=1$ it would say that the row data and the column data are
independent. Entropy duality turns the inequality into the first
estimate in \eqref{eq:bridge}.

To prove the correlation inequality, we decompose functions on the
board into irreducible representations of $S_N$, which are indexed by
the partitions of $N$, called shapes below, and bound the correlation
on each component. Shapes close to
a single row or column are handled by a cancellation involving marked
pairs of labels: an exchange can preserve the row data while changing
the sign of a matrix coefficient that depends only on the column
data. We bound the contribution of the exceptional configurations,
in which this exchange is unavailable. For the remaining shapes, we use their large
dimension and a bound on the fourth trace of the incidence matrix of
compatible row and column partitions. For odd $n$, the line groups
are alternating and the row and column data carry an additional
parity constraint within each line; the same argument keeps track of
these constraints. Section~\ref{glob:section} carries this out.

\paragraph{Comparison within a line.}
The second estimate in \eqref{eq:bridge} compares the complete
rearrangement of one line with the unit rotations of the shuffle. It
is proved by recursively bisecting each line. At an interval of
length $m$, the entropy chain rule separates the choice of labels in
its two halves from their internal arrangements, which are passed on
to the two halves. The first contribution is controlled by two
overlapping layers of short intervals of length of order $\sqrt m$,
each interval being rearranged at random. We compare these rearrangements with
parallel products of disjoint three-cycles, and build the three-cycles from
rotations. The triples in a layer have common spacings, so they can
all share the same row rotations; a further construction restores
every card outside the line.

At each node of the recursion, only columns meeting its interval
are used. A fixed column therefore enters only at the nodes on one
branch of the binary tree, where the interval lengths decrease
geometrically, so its costs form a geometric sum. Keeping this
dependence on the interval is what avoids an extra logarithmic
factor. Section~\ref{phys:section} gives the comparison, using a
smoothing estimate proved in Appendix~\ref{smooth:section}.

\paragraph{The matching bounds.}
The bounds in the other direction, $\alpha(P)\leq Cn^{-3}$ and
$\tmix^P(1/4)\geq cn^3\log n$, come from one function on the board.
Fix an initial board $x_0$, let $u_i(x)$ be the horizontal coordinate
of the card $i$ on the board $x$, and put
$F(x)=\sum_{i=1}^Na_i\cos(2\pi u_i(x)/n)$ with
$a_i=\cos(2\pi u_i(x_0)/n)$. Then $PF=(1-\gamma)F$ for a number
$\gamma\asymp n^{-3}$, and testing the logarithmic Sobolev inequality
on the functions $1+\varepsilon F$ gives $\alpha(P)\leq\gamma/2$. For
the mixing time we use the eigenfunction method of
Wilson~\cite{wilson}. By Lemma~\ref{lem:observable},
$F(x_0)\geq n^2/2$ and $F$ has mean zero under $\pi$, whereas the
standard deviation of $F$ is at most of order $n^{3/2}$, along the
chain and in equilibrium. The mean $(1-\gamma)^tF(x_0)$ along the
chain therefore exceeds these fluctuations during a constant multiple
of $\gamma^{-1}\log n\asymp n^3\log n$ updates, and until then the law
of the chain is far from $\pi$. The factor $\log n$ thus comes from
the ratio between $n^2$ and $n^{3/2}$.

\subsection{Organization and notation}\label{subsec:notation}
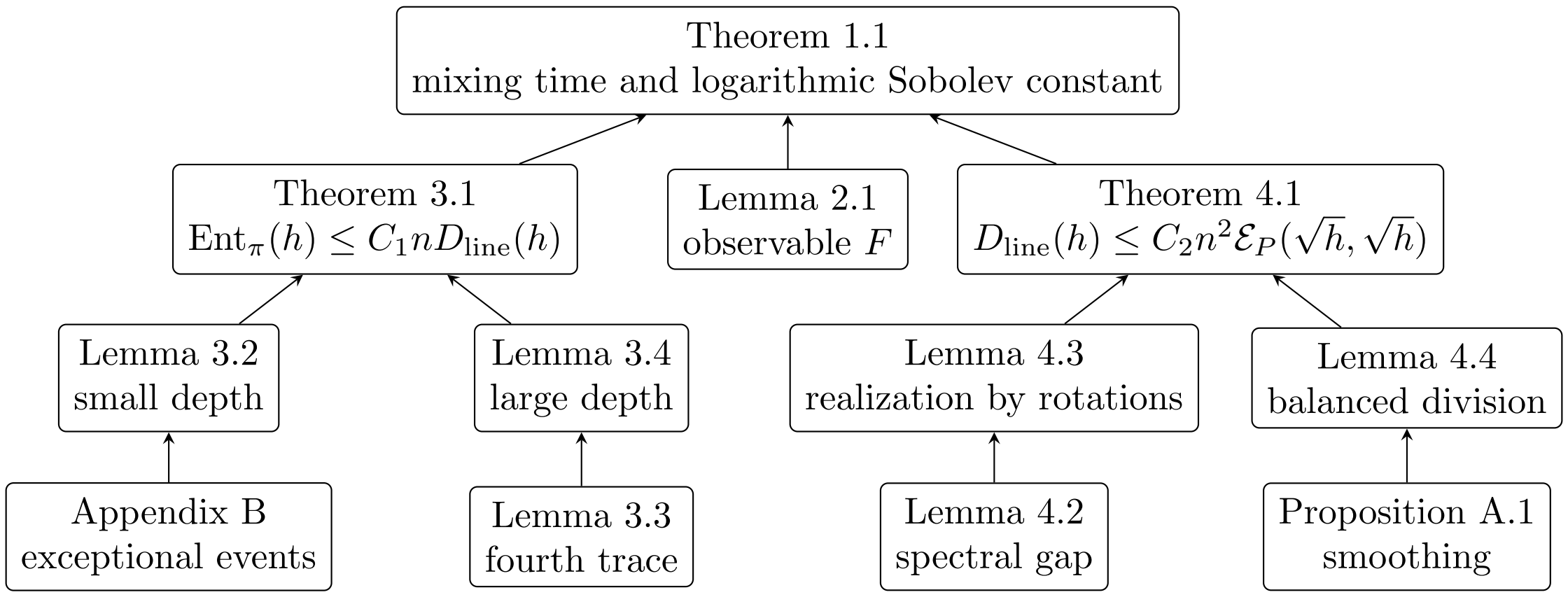
\begin{figure}[ht]
\centering
\begin{tikzpicture}[x=1cm,y=1cm,>=stealth,
  box/.style={draw,rounded corners=2pt,align=center,inner sep=4pt,
              font=\small}]
  \node[box] (main) at (0,4.5)
    {Theorem~\ref{thm:main}\\ mixing time and logarithmic Sobolev constant};
  \node[box] (glob) at (-3.9,3)
    {Theorem~\ref{thm:global}\\ $\Ent_\pi(h)\le C_1nD_{\rm line}(h)$};
  \node[box] (phys) at (3.9,3)
    {Theorem~\ref{thm:physical}\\
     $D_{\rm line}(h)\le C_2n^2\cE_P(\sqrt h,\sqrt h)$};
  \node[box] (obs) at (0,3)
    {Lemma~\ref{lem:observable}\\ observable $F$};
  \node[box] (low) at (-5.85,1.5) {Lemma~\ref{glob:low}\\ small depth};
  \node[box] (high) at (-1.95,1.5) {Lemma~\ref{glob:high}\\ large depth};
  \node[box] (rout) at (1.95,1.5)
    {Lemma~\ref{phys:routing}\\ realization by rotations};
  \node[box] (split) at (5.85,1.5)
    {Lemma~\ref{phys:split}\\ balanced division};
  \node[box] (exc) at (-5.85,0)
    {Appendix~\ref{exc:section}\\ exceptional events};
  \node[box] (trace) at (-1.95,0) {Lemma~\ref{glob:trace}\\ fourth trace};
  \node[box] (gap) at (1.95,0) {Lemma~\ref{phys:gap}\\ spectral gap};
  \node[box] (smooth) at (5.85,0)
    {Proposition~\ref{smooth:estimate}\\ smoothing};
  \draw[->] (glob) -- (main);
  \draw[->] (obs) -- (main);
  \draw[->] (phys) -- (main);
  \draw[->] (low) -- (glob);
  \draw[->] (high) -- (glob);
  \draw[->] (rout) -- (phys);
  \draw[->] (split) -- (phys);
  \draw[->] (exc) -- (low);
  \draw[->] (trace) -- (high);
  \draw[->] (gap) -- (rout);
  \draw[->] (smooth) -- (split);
\end{tikzpicture}
\caption{The main statements and their dependence. An arrow from one
statement to another means that the first is used in the proof of the
second.}
\label{fig:dependence}
\end{figure}

Subsection~\ref{subsec:classes} identifies the communicating classes.
Section~\ref{sec:mixing} deduces Theorem~\ref{thm:main} from
\eqref{eq:bridge} and proves the matching bounds, which rest on
Lemma~\ref{lem:observable}.
Section~\ref{glob:section} proves the first estimate in
\eqref{eq:bridge} and Section~\ref{phys:section} the second.
Appendix~\ref{smooth:section} contains the smoothing estimate used in
Section~\ref{phys:section}, and Appendix~\ref{exc:section} two
probability estimates used in the proof of Lemma~\ref{glob:low}.
Sections~\ref{glob:section}
and~\ref{phys:section} can be read independently of each other.
Figure~\ref{fig:dependence} shows how the main statements depend on
each other.

\paragraph{Conventions.}
A \emph{site} is an element of $V$. The site $(j,y)$ lies in column
$j$ and in row $y$, so that $j$ is its horizontal coordinate. We
identify a card with its label. A permutation $s$ of the sites acts on
boards by moving the cards, $(sx)(i)=s(x(i))$. We call such
permutations \emph{physical}, to distinguish them from permutations of
the labels. Products of permutations are composed from right to left:
in $st$ the permutation $t$ acts first. The letters $c$ and $C$ denote
positive absolute constants whose values may change from one
occurrence to the next; a subscript naming a parameter, as in
$C_{m_0}$, indicates dependence on that parameter. In
Section~\ref{glob:section},
where $C$ denotes the average over all columns, constants are written
$c$ or absorbed in the notation $O(\cdot)$. We write $a\asymp b$ if
$a\le Cb$ and $b\le Ca$. The notations $O(\cdot)$ and $o(\cdot)$ refer
to $n\to\infty$, with absolute implied constants.
Table~\ref{tab:notation} lists the main notation.

\begin{table}[ht]
\centering
\small
\begin{tabular}{@{}lp{0.60\textwidth}l@{}}
\hline
Symbol & Meaning & Defined in\\
\hline
$n$, $N$, $V$ & side of the board, number $n^2$ of cards, set of
  sites & \ref{subsec:model}\\
$G_n$ & group generated by the unit rotations of the lines &
  \ref{subsec:model}, \ref{subsec:classes}\\
$P$, $\pi$ & kernel of the half-lazy shuffle; uniform law on a
  communicating class & \ref{subsec:model}\\
$\E_\pi$, $\Ent_\pi$, $\cE_P$ & expectation and entropy under $\pi$;
  Dirichlet form & \ref{subsec:model}\\
$\E$, $\Ent$, $\|\cdot\|_2$ & without a subscript: taken under the law
  named at the start of the passage where they occur &
  \ref{glob:strategy}, \ref{phys:splitsection}, \ref{smooth:section}\\
$\alpha(P)$, $\tmix^P(\delta)$ & logarithmic Sobolev constant and
  mixing time & \ref{subsec:model}\\
$B_\ell$ & average over the line group of the line $\ell$ &
  \ref{subsec:outline}\\
$D_{\rm line}$ & mean entropy loss under the line averages &
  \ref{subsec:outline}\\
$u_i(x)$, $F$, $\gamma$ & horizontal coordinate of the card $i$ on the
  board $x$; Fourier observable and its rate of decay
  (Subsections~\ref{subsec:outline} and~\ref{subsec:lower} only) &
  \ref{subsec:lower}\\
$R$, $C$ & averages over all rows and over all columns
  (Section~\ref{glob:section} only) & \ref{glob:setting}\\
$\lambda$, $d_\lambda$, $k$, $\beta_\lambda$ & shape, dimension of
  its Specht module, depth, correlation of $R$ and $C$ on the
  $\lambda$-component & \ref{glob:strategy}\\
$U_s$, $e_s(f)$ & action of a site permutation $s$ on functions, and
  its energy & \ref{phys:section}\\
$R_y$, $C_j$ & unit rotations of the row $y$ and of the column $j$
  (Section~\ref{phys:section} only) & \ref{phys:section}\\
$m$, $b$ & length of an interval of a line; length scale of the short
  intervals & \ref{phys:section}\\
$B_{\mathcal I}^+$, $K_{\mathcal I}$ & average over the alternating
  groups of the short intervals in $\mathcal I$; random update by
  disjoint three-cycles & \ref{phys:gapsection}\\
\hline
\end{tabular}
\caption{Main notation.}\label{tab:notation}
\end{table}

\subsection{Communicating classes}\label{subsec:classes}
We identify the group $G_n$ generated by the unit rotations. The
answer is known from work on the torus puzzle \cite{amano}, as noted
in \cite[Section~6.1]{bms}. We include a proof. It uses the following
elementary fact, which is applied again
to the line groups in Section~\ref{glob:section}.

\begin{lemma}\label{lem:triples}
Let $W$ be a finite set and let $\mathcal T$ be a family of
three-element subsets of $W$ whose union is $W$. Suppose that
$\mathcal T$ is connected: any two members of $\mathcal T$ are joined
by a chain of members of $\mathcal T$ in which consecutive members
intersect. Then the three-cycles supported on the members of
$\mathcal T$ generate the alternating group on $W$.
\end{lemma}
\begin{proof}
Write $A_S$ for the alternating group on a subset $S$ of $W$. Start
with a member $S$ of $\mathcal T$; its two three-cycles generate
$A_S$. Suppose now that $A_S$ is already available for some
$S\subseteq W$ that is a union of members of $\mathcal T$. If
$S\ne W$, connectedness gives a member $T$ of $\mathcal T$ that meets
$S$ and is not contained in $S$.

If $T$ meets $S$ in two sites and $x$ is its third site, we can
adjoin $x$. Indeed, $A_S$ is transitive on unordered pairs of $S$, so
conjugating the three-cycles on $T$ by $A_S$ gives a three-cycle on
$\{a,b,x\}$ for every pair $a,b$ in $S$, and inverses supply either
orientation. Every three-cycle on $S\cup\{x\}$ lies in $A_S$ or has
this form, so $A_{S\cup\{x\}}$ is available.

If $T=\{a,x,y\}$ meets $S$ only at $a$, conjugate the cycle
$(a\,x\,y)$ by an element of $A_S$ to obtain $(b\,x\,y)$ with
$b\ne a$ in $S$. Multiplying the first cycle
by the inverse of the second gives a three-cycle on $\{a,b,x\}$, which
meets $S$ in two sites. Now adjoin $x$ as before and then $y$.
Connectedness allows this procedure to reach every element of $W$.
\end{proof}

\begin{proposition}\label{prop:classes}
For every $n\ge2$, the group $G_n$ is the symmetric group $S_N$ on
$V$ if $n$ is even, and the alternating group $A_N$ if $n$ is odd.
\end{proposition}
\begin{proof}
Take a nonidentity power of the unit rotation of a row and a
nonidentity power of the unit rotation of a crossing column. Their
supports meet at one site, the corner, and their commutator cycles
the corner and one further site of each of the two lines. Varying the powers gives
every such corner triple, in either orientation. The resulting family
of triples covers $V$ and is connected, so $G_n$ contains $A_N$ by
Lemma~\ref{lem:triples}.

A unit rotation of a line has sign $(-1)^{n-1}$, so every generator
is even for odd $n$, whereas an odd generator is available for even
$n$. This proves the assertion.
\end{proof}

\section{Entropy decay and mixing}\label{sec:mixing}
This section proves Theorem~\ref{thm:main} in two parts.
Subsection~\ref{subsec:decay} uses the two estimates in
\eqref{eq:bridge}, which are proved in Sections~\ref{glob:section}
and~\ref{phys:section}. The estimates give a logarithmic Sobolev
inequality, hence exponential decay of entropy and the upper bound on
the mixing time. Subsection~\ref{subsec:lower} proves the bounds in the
other direction from the Fourier observable of
Lemma~\ref{lem:observable}, independently of \eqref{eq:bridge}.

\subsection{Entropy decay}\label{subsec:decay}
We first prove the bound $\alpha(P)\geq cn^{-3}$ and the upper bounds
on the mixing time.

\begin{proof}[Proof of Theorem~\ref{thm:main}, first part]
\step{Step 1: Logarithmic Sobolev inequality.}
Let $f$ be a real function and put $h=f^2$, so that $\sqrt h=|f|$.
Since $\cE_P(|f|,|f|)\leq\cE_P(f,f)$, the two estimates in
\eqref{eq:bridge} give
\begin{equation}\label{eq:lsi}
 \Ent_\pi(f^2)\leq C n^3\cE_P(f,f).
\end{equation}
This proves $\alpha(P)\geq cn^{-3}$.

\step{Step 2: Entropy contraction in discrete time.}
To apply \eqref{eq:lsi} in discrete time we follow an argument of
Miclo \cite{miclo}. For every
probability measure $\mu$ and every $h\geq0$, we have
$\Ent_\mu(h)\geq\Var_\mu(\sqrt h)$, where $\Ent_\mu$ and $\Var_\mu$
are the entropy and the variance under $\mu$. After normalizing $h$ to
have mean one under $\mu$, this follows by integrating
$z\log z-z+1\geq(\sqrt z-1)^2$; the case of the zero function is
immediate.
By stationarity and reversibility, for a density $h$,
\begin{align*}
 \Ent_\pi(h)-\Ent_\pi(Ph)
 &=\sum_x\pi(x)\Ent_{P(x,\cdot)}(h)\\
 &\geq\sum_x\pi(x)\Var_{P(x,\cdot)}(\sqrt h)\\
 &=\langle\sqrt h,(I-P^2)\sqrt h\rangle_\pi
 \geq \cE_P(\sqrt h,\sqrt h).
\end{align*}
The first line holds because $\pi$ is stationary: $P(h\log h)$ has the
same mean under $\pi$ as $h\log h$, and $Ph$ has the same mean as $h$.
The second line is the inequality above with $\mu=P(x,\cdot)$. In the
third line, stationarity gives
$\sum_x\pi(x)\Var_{P(x,\cdot)}(\sqrt h)
=\|\sqrt h\|_{2,\pi}^2-\|P\sqrt h\|_{2,\pi}^2$, and reversibility
turns the last norm into $\langle\sqrt h,P^2\sqrt h\rangle_\pi$.
The last inequality uses the half-laziness, which puts the spectrum of
$P$ in $[0,1]$, so that $I-P^2=(I-P)(I+P)$ dominates $I-P$. Combining
this one-step loss with \eqref{eq:lsi},
applied to $f=\sqrt h$, shows that entropy contracts by a factor at
most $1-c/n^3$ at each update:
\[
 \Ent_\pi(Ph)\leq\Ent_\pi(h)-\cE_P(\sqrt h,\sqrt h)
 \leq(1-c/n^3)\Ent_\pi(h).
\]

\step{Step 3: Upper bound on the mixing time.}
Let $x$ be a board, let $\pi$ be the uniform law on its class, and
let $h_0$ be the density of the point mass at $x$ with respect to
$\pi$. By reversibility, $P^th_0$ is the density of $P^t(x,\cdot)$.
The density $h_0$ has entropy
$\log|G_n|\leq 2n^2\log n$. By Step~2,
\[
 \Ent_\pi(P^t h_0)\leq e^{-ct/n^3}\log|G_n|.
\]
Pinsker's inequality gives
$\|P^t(x,\cdot)-\pi\|_{\rm TV}^2
 \leq \frac12 e^{-ct/n^3}\log|G_n|$.
Taking $t\geq Cn^3(\log n+\log(1/\delta))$ makes the right-hand
side at most $\delta^2$, uniformly in the initial board, and proves
the asserted upper bounds.
\end{proof}

\subsection{The matching bounds}\label{subsec:lower}
It remains to prove that $\alpha(P)\leq Cn^{-3}$ and that
$\tmix^P(1/4)\geq cn^3\log n$. For the mixing time we follow the
eigenfunction method of Wilson~\cite{wilson}. The method uses an
eigenfunction of $P$ with three properties: its value at the initial
board is large, its eigenvalue is close to one, and its fluctuations
along the chain are small.

Fix a board $x_0$ and let $\pi$ be the uniform law on its class. For a
board $x$ of this class and a label $i$, let
$u_i(x)\in\mathbb Z/n\mathbb Z$ be the horizontal coordinate of the
card $i$ on $x$. Put
\begin{equation}\label{eq:observable}
 a_i=\cos\frac{2\pi u_i(x_0)}{n},\qquad
 F(x)=\sum_{i=1}^N a_i\cos\frac{2\pi u_i(x)}{n},\qquad
 \gamma=\frac{1-\cos(2\pi/n)}{4n}.
\end{equation}
The weights $a_i$ make every card contribute to $F(x_0)$ with the same
sign. With equal weights the sum would be zero on every board, because
each column always holds $n$ cards. Since
$2\varphi^2/\pi^2\leq1-\cos\varphi\leq\varphi^2/2$ for
$0\leq\varphi\leq\pi$, the choice $\varphi=2\pi/n$ gives
\begin{equation}\label{eq:gamma}
 \frac{2}{n^3}\leq\gamma\leq\frac{\pi^2}{2n^3}
 \qquad\text{for every }n\geq2.
\end{equation}

\begin{lemma}\label{lem:observable}
Let $n\geq2$, and let $F$ and $\gamma$ be defined by
\eqref{eq:observable}.
\begin{enumerate}
\item[(i)] $F(x_0)\geq n^2/2$, with equality for $n\geq3$, and
$\E_\pi[F]=0$.
\item[(ii)] $PF=(1-\gamma)F$.
\item[(iii)] Let $(X_t)_{t\geq0}$ be the chain with kernel $P$ started
at $x_0$, or started from $\pi$. Then
$\Var(F(X_t))\leq4\pi^2/\gamma$ for every $t\geq0$.
\end{enumerate}
\end{lemma}
\begin{proof}
\step{Part (i).}
At $x_0$ every column holds $n$ cards, and a card $i$ in column $j$
has $a_i=\cos(2\pi j/n)$ and contributes $a_i^2$ to $F(x_0)$. Hence
\[
 F(x_0)=n\sum_{j=0}^{n-1}\cos^2\frac{2\pi j}{n}
 =\frac{n^2}{2}+\frac n2\sum_{j=0}^{n-1}\cos\frac{4\pi j}{n}.
\]
The last sum is zero for $n\geq3$ and equals $2$ for $n=2$. The
second assertion is proved in part (ii).

\step{Part (ii).}
Fix a board $x$ and a label $i$, and write $\theta=2\pi u_i(x)/n$. One
update moves the card $i$ one step to the right with probability
$1/(8n)$ and one step to the left with the same probability, namely
when its row is rotated, and it leaves $u_i$ unchanged otherwise.
Since
$\cos(\theta+2\pi/n)+\cos(\theta-2\pi/n)=2\cos(2\pi/n)\cos\theta$,
the mean of $\cos(2\pi u_i/n)$ after one update from $x$ is
\[
 \Bigl(1-\frac1{4n}\Bigr)\cos\theta
 +\frac1{4n}\cos\frac{2\pi}{n}\cos\theta
 =(1-\gamma)\cos\theta .
\]
For $n=2$ the two steps lead to the same column, and the display is
unchanged. Summing over $i$ with the weights $a_i$ gives
$PF=(1-\gamma)F$. Since $\pi$ is stationary,
$\E_\pi[F]=\E_\pi[PF]=(1-\gamma)\E_\pi[F]$, and $\gamma>0$ gives
$\E_\pi[F]=0$.

\step{Part (iii).}
Put $\sigma_t^2=\Var(F(X_t))$. By part (ii),
$\E[F(X_{t+1})\mid X_t]=(1-\gamma)F(X_t)$. The decomposition of the
variance according to $X_t$ therefore reads
\[
 \sigma_{t+1}^2
 =\Var\bigl(\E[F(X_{t+1})\mid X_t]\bigr)
  +\E\bigl[\Var(F(X_{t+1})\mid X_t)\bigr]
 =(1-\gamma)^2\sigma_t^2+\E\bigl[\Var(F(X_{t+1})\mid X_t)\bigr].
\]
The conditional variance is at most the conditional mean of
$(F(X_{t+1})-F(X_t))^2$, and one update changes $F$ by at most $2\pi$.
Indeed, a column rotation changes no horizontal coordinate, and a row
rotation moves $n$ cards by one step, so that $n$ terms of $F$ change
by at most $2\pi/n$ each, because $|a_i|\leq1$ and the cosine is
$1$-Lipschitz. Hence
$\sigma_{t+1}^2\leq(1-\gamma)^2\sigma_t^2+4\pi^2$. For the chain
started at $x_0$ we have $\sigma_0^2=0$, and iterating gives
\[
 \sigma_t^2\leq4\pi^2\sum_{s\geq0}(1-\gamma)^{2s}
 =\frac{4\pi^2}{\gamma(2-\gamma)}\leq\frac{4\pi^2}{\gamma}.
\]
For the chain started from $\pi$, the variance $\sigma_t^2=\Var_\pi(F)$
does not depend on $t$, and the same inequality gives
$\gamma(2-\gamma)\Var_\pi(F)\leq4\pi^2$, hence the same bound.
\end{proof}

\begin{proof}[Proof of Theorem~\ref{thm:main}, second part]
\step{Step 4: Upper bound on the logarithmic Sobolev constant.}
In the variational formula for $\alpha(P)$, take $f=1+\varepsilon F$
with $\varepsilon>0$. By Lemma~\ref{lem:observable}, $F$ has mean zero
under $\pi$, it is not constant because $F(x_0)>0$, and
$(I-P)F=\gamma F$. Hence
$\cE_P(f,f)=\varepsilon^2\langle F,(I-P)F\rangle_\pi
=\varepsilon^2\gamma\Var_\pi(F)$. We expand the entropy as
$\varepsilon\to0$ with $n$ fixed. Since
$\log f^2=2\varepsilon F-\varepsilon^2F^2+O(\varepsilon^3)$ uniformly
on the class,
\[
 f^2\log f^2=2\varepsilon F+3\varepsilon^2F^2+O(\varepsilon^3),
 \qquad
 \E_\pi[f^2]\log\E_\pi[f^2]
 =\varepsilon^2\E_\pi[F^2]+O(\varepsilon^4),
\]
so that $\Ent_\pi(f^2)=2\varepsilon^2\Var_\pi(F)+O(\varepsilon^3)$.
The ratio $\cE_P(f,f)/\Ent_\pi(f^2)$ tends to $\gamma/2$. This gives,
with \eqref{eq:gamma},
\[
 \alpha(P)\leq\frac\gamma2\leq\frac{\pi^2}{4n^3}
 \qquad\text{for every }n\geq2.
\]

\step{Step 5: Separation from equilibrium.}
Let $(X_t)$ be the chain started at $x_0$, and put
\[
 t=\Bigl\lfloor\frac{\log n}{8\gamma}\Bigr\rfloor,\qquad
 M=\E[F(X_t)].
\]
By parts (i) and (ii) of Lemma~\ref{lem:observable},
$M=(1-\gamma)^tF(x_0)\geq(1-\gamma)^tn^2/2$. The definition of
$\gamma$ gives $\gamma\leq1/(2n)\leq1/4$, so that
$\log(1-\gamma)\geq-2\gamma$ and
$(1-\gamma)^t\geq e^{-2\gamma t}\geq n^{-1/4}$. Hence
\begin{equation}\label{eq:mean}
 M\geq\tfrac12n^{7/4}.
\end{equation}
Let $A$ be the set of boards $x$ with $F(x)>M/2$. Chebyshev's
inequality and part (iii) of the lemma give, for the chain started at
$x_0$,
\[
 P^t(x_0,A^c)\leq\mathbb P\bigl(|F(X_t)-M|\geq M/2\bigr)
 \leq\frac{4\Var(F(X_t))}{M^2}\leq\frac{16\pi^2}{\gamma M^2}.
\]
Since $F$ has mean zero under $\pi$, the same argument gives
$\pi(A)\leq16\pi^2/(\gamma M^2)$. By \eqref{eq:gamma} and
\eqref{eq:mean},
\[
 \|P^t(x_0,\cdot)-\pi\|_{\rm TV}\geq P^t(x_0,A)-\pi(A)
 \geq1-\frac{32\pi^2}{\gamma M^2}\geq1-\frac{64\pi^2}{\sqrt n}.
\]
Since $64\pi^2<750$, the right-hand side is larger than $1/4$ for
$n\geq10^6$. The distance to equilibrium does not increase with time,
so $\tmix^P(1/4)>t$ for these $n$. By \eqref{eq:gamma} again,
\[
 \tmix^P(1/4)>\frac{\log n}{8\gamma}-1\geq\frac{n^3\log n}{4\pi^2}-1
 \qquad\text{for every }n\geq10^6.
\]
For $2\leq n<10^6$ the distance to equilibrium at time zero is
$1-1/|G_n|\geq1/2$, so that $\tmix^P(1/4)\geq1$. The lower bound of
the theorem holds for these finitely many $n$ after decreasing the
constant $c$.
\end{proof}

\begin{remark}[Two-step chain]
The proof of the upper bound in \cite[Section~8]{bms} works with the
chain that performs
two updates at a time, whose kernel is $P^2$. The mixing-time bounds
of Theorem~\ref{thm:main} hold for that chain with the same orders:
total variation distance to stationarity is nonincreasing, so for every
$\delta$ the exact time conversion is
$\tmix^{P^2}(\delta)=\lceil\tmix^P(\delta)/2\rceil$.
\end{remark}

\section{Entropy factorization over rows and columns}\label{glob:section}

This section proves the first estimate in \eqref{eq:bridge}. When the
line groups are symmetric, rearranging every row at random retains
only the set of labels in each row, and rearranging every column
retains only the set of labels in each column. We call these the row
data and the column data of the board. A single label has independent
row and column coordinates under the uniform board law. The full row
data and column data are dependent, however, since every row set must
meet every column set in exactly one label. We bound this dependence
to obtain an entropy factorization uniform in the board size.

\subsection{Setting and statement}\label{glob:setting}

Throughout this section the line groups are either all symmetric,
with $n\ge2$, or all alternating, with $n\ge3$. The torus shuffle has
symmetric line groups for even $n$ and alternating line groups for
odd $n$, but the parity of $n$ plays no role in this section.
Write $B_\ell$ for uniform averaging over the line group of the
line $\ell$. The law $\pi$ is uniform on a single class, that is, on
an orbit of the group generated by the line groups.
The symmetric choice has one class, and the alternating choice has the two
global parity classes. Indeed, site transpositions on a connected graph
generate the symmetric group. For the alternating case take all triples
contained in a row or column. The resulting family of triples is connected,
so its three-cycles generate the alternating group, by
Lemma~\ref{lem:triples}.

The parallel row and column averages are the orthogonal projections
\[
 R=\prod_{\ell\text{ a row}}B_\ell,\qquad C=\prod_{\ell\text{ a column}}B_\ell .
\]
Each product averages independently within the lines of its orientation.
Thus $R$ and $C$ are conditional expectations given the row orbit and
the column orbit of the board, that is, its orbit under the product
of the row groups and under the product of the column groups. For
symmetric line groups these orbits are determined by the row data and
by the column data. For alternating line groups they also
retain within-line parity information. A function $F$ satisfies
$F=RF$ exactly when it depends on the board only through its row
orbit, and $H=CH$ exactly when $H$ depends only on the column orbit.

\begin{theorem}\label{thm:global}
There is an absolute $p\in(1,2)$ such that for every nonnegative
$F=RF$ and $H=CH$,
\begin{equation}\label{glob:bilinear}
 \E_\pi[FH]\le \|F\|_{L^p(\pi)}\|H\|_{L^p(\pi)}.
\end{equation}
In particular, for every $h\ge0$,
\begin{equation}\label{glob:factorization}
 (2-p)\Ent_\pi(h)
 \le 2\Ent_\pi(h)-\Ent_\pi(Rh)-\Ent_\pi(Ch)
 \le 2n D_{\rm line}(h),
\end{equation}
where $D_{\rm line}(h)=(2n)^{-1}\sum_\ell
[\Ent_\pi(h)-\Ent_\pi(B_\ell h)]$.
\end{theorem}

The outer terms of \eqref{glob:factorization} give the first estimate
in \eqref{eq:bridge} with $C_1=2/(2-p)$. The proof deduces
\eqref{glob:factorization} from \eqref{glob:bilinear} in its last
step, by entropy duality and tensorization.

The first inequality in \eqref{glob:factorization} says that
$\Ent_\pi(Rh)+\Ent_\pi(Ch)\leq p\Ent_\pi(h)$. For $p=2$ this is
trivial, since averaging does not increase entropy. For $p=1$ it would
hold if the row orbit and the column orbit were independent, and it
does hold when $h$ depends on the site of a single card. The theorem
asserts that some $p<2$ works for every $n$, although the two orbits
are dependent. Applied to $Rh$, the inequality gives
$\Ent_\pi(CRh)\leq(p-1)\Ent_\pi(Rh)$: once the rows have been
rearranged, rearranging the columns removes a fixed fraction of the
remaining entropy.

\begin{remark}[The factor $n$]
The factor $n$ in \eqref{glob:factorization} has the correct order.
Let $h$ be $n^2$ times the indicator
that a fixed card occupies a prescribed site. Then
$\Ent_\pi(h)=2\log n$.
Only the row and column of that site contribute to $D_{\rm line}$,
and each of the two line averages decreases the entropy by $\log n$.
Hence $D_{\rm line}(h)=(\log n)/n$ and
$\Ent_\pi(h)=2nD_{\rm line}(h)$.
\end{remark}

\subsection{Strategy and notation}\label{glob:strategy}

We first work on all of $S_N$, where $N=n^2$, with the uniform probability,
including when the line groups are alternating; the passage to a
single class is made at the end of the proof. Until that passage,
$\E$ and $\mathbb P$, and the norms and inner products of functions,
refer to this uniform probability unless another law is named; after
it they refer to $\pi$. An element $g$ of $S_N$
is a board: it maps labels to sites. Physical permutations act on
the left and label permutations on the right. We use the real Specht
representations
$(V_\lambda,\rho_\lambda)$ for $\lambda\vdash N$, the Young branching
rule and the hook-length formula as in \cite{etingof}.
Let $d_\lambda=\dim V_\lambda$ and let $\lambda'$ denote the transpose
partition. A partition of $N$ is also called a shape. The \emph{depth}
of the shape $\lambda$ is
\[
 k=N-\max(\lambda_1,\lambda'_1),
\]
the number of boxes outside the longer of the first row and the first
column. The operators $R$ and $C$ are averages of left translations,
so they preserve the $\lambda$-isotypic component of $L^2(S_N)$.
On this component write $\beta_\lambda=\|RC\|_{2\to2}$.
For two orthogonal projections this norm is the largest absolute inner
product between unit vectors in their respective ranges. We set it
equal to zero if one range is zero. Thus $\beta_\lambda$ measures the
correlation between the row data and the column data on one
irreducible component.

To pass from $L^2$
correlation bounds to the $L^p$ estimate in \eqref{glob:bilinear},
we use interpolation with a dimension factor. For
$\theta=2/p-1$, \eqref{glob:interp} reduces the problem to proving
\[
 \sum_{\lambda:\,k\ge1}d_\lambda^{2\theta}\beta_\lambda=o(1)
\]
for a fixed $p<2$ as $n$ tends to infinity. The sum excludes the constant component
and, for alternating line groups, the global parity component. We
treat two regimes of the depth, separated by the threshold
$n^{21/20}$. Its exponent lies between
two constraints: the proof of Lemma~\ref{glob:high} needs the depth to
be large compared with $n\log^2n$, and the proof of
Lemma~\ref{glob:low} needs it to be small compared with $n^{5/4}$
(Appendix~\ref{exc:section}).

For shapes of small depth, $1\le k\le n^{21/20}$, we represent a
matrix coefficient using at most $2k$ marked labels.
Exchanging the column coordinates of a suitable pair preserves the row
data and reverses the sign on the column side
(Figure~\ref{fig:exchange}). We bound the
exceptional events conditionally on the row data or on the column
data. The resulting bound is Lemma~\ref{glob:low}.

For shapes of large depth, $k>n^{21/20}$, we use the multiplicity
$d_\lambda$ in the regular
representation. A bound on the total fourth trace then controls
$\beta_\lambda$ in terms of $d_\lambda$.
Lemmas~\ref{glob:trace} and~\ref{glob:high} carry this out.

Subsection~\ref{glob:interpolation} combines these bounds, handles the
finitely many remaining board sizes, and passes from centered correlation
to \eqref{glob:bilinear}. Entropy duality then gives
\eqref{glob:factorization}.

\subsection{Representations of small depth}\label{glob:lowsection}

\begin{lemma}\label{glob:low}
For all sufficiently large $n$ and $1\le k\le n^{21/20}$,
$\beta_\lambda\le n^{-k/1000}$.
\end{lemma}

\begin{proof}
Fix a shape $\lambda$ of depth $k$ with $1\le k\le n^{21/20}$. When an
estimate below is stated for all sufficiently large $n$, the threshold
is an absolute constant, independent of $k$ and of $\lambda$. If
$\beta_\lambda=0$ there is nothing to prove, so we assume that the
ranges of $R$ and of $C$ on the $\lambda$-isotypic component are both
nonzero.

Steps~1 to~3 reduce the lemma to a bound on $\E[FH]$ for two
functions $F$ and $H$ of the sites of at most $2k$ labels, which we
call marked. Step~5 defines an exchange of the columns of two marked
labels that keeps $F$ and reverses the sign of $H$; see
Figure~\ref{fig:exchange}. The configurations of the marked labels in
which this exchange is possible therefore contribute nothing to
$\E[FH]$. Steps~4 and~6 show that the other configurations are rare,
even when the rows of the marked labels, or their columns, are given.
Step~7 concludes by the Cauchy--Schwarz inequality.

\step{Step 1: Reduction to a long first row.}
We may assume that $\lambda_1=N-k$, that is,
$\lambda_1\ge\lambda'_1$. The $k$ boxes counted by the depth then lie
below the first row.

For alternating line groups, multiplying a function on $S_N$ by the
sign character is an isometry of $L^2(S_N)$. It maps the
$\lambda$-isotypic component onto the $\lambda'$-isotypic component,
because tensoring a Specht module with the sign character transposes
its shape. It commutes with $R$ and with $C$, because every element
of a line group is an even permutation. Hence
$\beta_{\lambda'}=\beta_\lambda$. Since $\lambda'$ has the same depth,
we may replace $\lambda$ by $\lambda'$ if $\lambda_1<\lambda'_1$.

For symmetric line groups, a function in the range of $R$ depends on
the board only through the row of each label, that is, through a
word of length $N$ in $n$ letters. The functions of such words form
the $N$-fold tensor power of an $n$-dimensional space. Label
permutations permute its factors, and the antisymmetrizer of any
$n+1$ factors vanishes on it. A Specht module whose shape has more
than $n$ rows is generated by a polytabloid that is antisymmetric in
the labels of its first column, which number more than $n$. Every
homomorphism from this module into the tensor power is therefore
zero. The isotypic components of $L^2(S_N)$ are the same for physical
permutations and for label permutations, and the range of $R$ meets
the $\lambda$-isotypic component by assumption. Hence
$\lambda'_1\le n$. Since $N-k\ge n^2-n^{21/20}>n$ for all
sufficiently large $n$, the maximum in the definition of $k$ is
attained by $\lambda_1$.

\step{Step 2: Marked labels.}
Fix a tableau of shape $\lambda$, as in the left panel of
Figure~\ref{fig:exchange}. A label is \emph{marked} if its column in
the tableau has height at least two. In each such column we pair as
many labels as possible. These pairs are the \emph{designated pairs},
and from now on a pair always means a designated pair. A column of odd
height leaves one marked label unpaired, which we call a
\emph{singleton}. The pairs and the singletons are the
\emph{clusters}, so every marked label belongs to exactly one cluster.

Let $t$ be the number of marked labels and $s$ the number of pairs.
A column of height $a\ge2$ contributes $a-1$ boxes below the first row,
$a$ marked labels, and $\lfloor a/2\rfloor$ pairs, and these numbers
satisfy
\[
 a-1\le a\le2(a-1),\qquad (a-1)/2\le\lfloor a/2\rfloor\le a-1 .
\]
By Step~1 there are $k$ boxes below the first row, so summing over
the columns of height at least two gives
\begin{equation}\label{glob:counts}
 k\le t\le2k,\qquad k/2\le s\le k .
\end{equation}
In the figure, $k=3$, $t=5$ and $s=2$.

Recall that the polytabloid of the tableau is obtained by applying
its column antisymmetrizer to its row tabloid. Let $v\in V_\lambda$ be
this polytabloid, normalized to unit length. We use two properties of
$v$. First, $v$ changes sign
under the transposition of the two labels of a pair, because these
labels lie in the same column of the tableau. Second, every
permutation of the unmarked labels fixes $v$. Indeed, the unmarked
labels lie in the first row of the tableau, so such a permutation
fixes the row tabloid, and it commutes with the column
antisymmetrizer, which moves only marked labels.

\step{Step 3: Matrix coefficients.}
For $w\in V_\lambda$ consider the matrix coefficient
\[
 f_w(g)=\sqrt{d_\lambda}\,\langle w,\rho_\lambda(g)v\rangle .
\]
For a physical permutation~$\sigma$ and a label permutation~$\tau$,
\[
 f_w(\sigma^{-1}g)=f_{\rho_\lambda(\sigma)w}(g),\qquad
 f_w(g\tau)=\sqrt{d_\lambda}\,
 \langle w,\rho_\lambda(g)\rho_\lambda(\tau)v\rangle .
\]
The first identity holds because $\rho_\lambda(\sigma)$ is orthogonal.
Thus physical permutations act on the vector $w$ and label
permutations on $v$.

Let~$P_{\rm row}$ and~$P_{\rm col}$ be the averages
of~$\rho_\lambda(\sigma)$ over the product of the row groups and
over the product of the column groups, respectively. Averaging the
first identity over~$\sigma$ yields
\begin{equation*}
 Rf_w=f_{P_{\rm row}w},\qquad Cf_w=f_{P_{\rm col}w}\,.
\end{equation*}
By Schur orthogonality, $\E[f_wf_{w'}]=\langle w,w'\rangle$ for all
$w,w'\in V_\lambda$, so the map~$w\mapsto f_w$ is an isometry onto a
copy of~$V_\lambda$ inside~$L^2(S_N)$. The copies obtained in the same
way from the vectors of an orthonormal basis of~$V_\lambda$ that
contains~$v$ are mutually orthogonal and span the~$\lambda$-isotypic
component. On each of them $R$ and $C$ act through~$P_{\rm row}$
and~$P_{\rm col}$.
Hence~$\beta_\lambda=\|P_{\rm row}P_{\rm col}\|_{2\to2}$ is the
largest value of $|\E[FH]|$ over the coefficients $F=f_w$ and
$H=f_{w'}$ of unit norm that satisfy $F=RF$ and $H=CH$. It therefore
suffices to prove that
\[
 |\E[FH]|\le n^{-k/1000}
\]
for all such $F$ and $H$ and all sufficiently large $n$. We fix $F$
and $H$ for the rest of the proof.

We now use the two properties of $v$. The \emph{configuration} of a
board $g$ is its restriction $\eta$ to the marked labels, so
$\eta$ sends each marked label to its site. Two boards with the
same configuration differ by a permutation of the unmarked labels,
which fixes $v$. By the second identity above, $F$ and $H$ therefore
depend on $g$ only through $\eta$, and we write $F(\eta)$ and
$H(\eta)$. If $\tau$ is the transposition of the two labels of a
pair, then the configuration of $g\tau$ is $\eta\circ\tau$ and
$\rho_\lambda(\tau)v=-v$, so that
\[
 F(\eta\circ\tau)=-F(\eta),\qquad H(\eta\circ\tau)=-H(\eta).
\]
Under the uniform law on $S_N$, the configuration is a uniformly
distributed injection of the set of marked labels into the set of
sites.

\step{Step 4: Conditional laws of the marked labels.}
The functions $F$ and $H$ have unit norm in $L^2$ but no uniform
bound, so the probability of an event does not control its
contribution to $\E[FH]$. We instead bound conditional
probabilities given the following orbits.

The \emph{row orbit} of a configuration $\eta$ is its orbit under the
product of the row groups, where a physical permutation $\sigma$ maps
$\eta$ to $\sigma\circ\eta$. The \emph{column orbit} is defined in the
same way. Since $F=RF$, the function $F$ is invariant under the row
groups, so $F(\eta)$ depends only on the row orbit of $\eta$.
Likewise $H(\eta)$ depends only on its column orbit. Given its row
orbit, $\eta$ is uniformly distributed on this orbit.

For symmetric line groups, two configurations lie in the same row
orbit exactly when they assign the same row to every marked label.
Given the row orbit, the marked labels of each row thus occupy
uniformly chosen distinct columns, and different rows are
independent. We call this law, which depends only on the rows of the
marked labels, the \emph{symmetric conditional law}.

For alternating line groups, call a site \emph{unmarked} if it
carries no marked label, and consider two configurations that assign
the same row to every marked label. In each row some permutation of
the sites carries the first configuration to the second. If the row
has at least two unmarked sites, this permutation can be chosen even,
because composing it with the transposition of two unmarked sites
changes its parity and moves no marked label. Hence the two
configurations lie in the same row orbit if they differ only on rows
with at least two unmarked sites. In a row with $n-1$ or $n$ marked
labels the permutation is unique, and the ways of placing these
labels in the row form two orbits of the alternating group of the
row, of equal size.

Given a row orbit, different rows are therefore independent. In a row
with at least two unmarked sites the marked labels occupy uniformly
chosen distinct columns, as in the symmetric case. In each of the at
most $t/(n-1)$ other rows they are uniformly distributed on one of
the two orbits. Every event thus has conditional probability at most
$2^{t/(n-1)}$ times its probability under the symmetric conditional
law for the same rows. We call this the density bound. The same
statements hold with rows and columns interchanged.

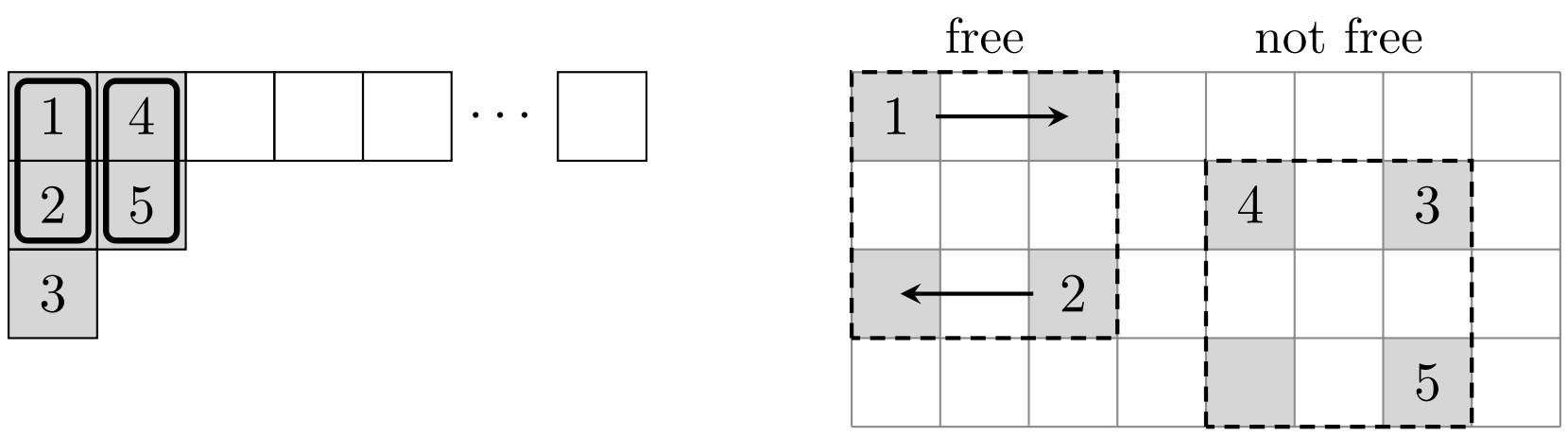
\begin{figure}[ht]
\centering
\begin{tikzpicture}[x=0.62cm,y=0.62cm,>=stealth]
  \begin{scope}[yshift=2.48cm]
    \fill[black!16] (0,0) rectangle (1,-3);
    \fill[black!16] (1,0) rectangle (2,-2);
    \foreach \x in {0,...,4} \draw (\x,0) rectangle (\x+1,-1);
    \node at (5.6,-0.5) {$\cdots$};
    \draw (6.2,0) rectangle (7.2,-1);
    \draw (0,-1) rectangle (1,-2);
    \draw (1,-1) rectangle (2,-2);
    \draw (0,-2) rectangle (1,-3);
    \node at (0.5,-0.5) {$1$};
    \node at (0.5,-1.5) {$2$};
    \node at (0.5,-2.5) {$3$};
    \node at (1.5,-0.5) {$4$};
    \node at (1.5,-1.5) {$5$};
    \draw[very thick,rounded corners=2pt] (0.1,-0.1) rectangle (0.9,-1.9);
    \draw[very thick,rounded corners=2pt] (1.1,-0.1) rectangle (1.9,-1.9);
  \end{scope}
  \begin{scope}[xshift=5.9cm]
    \foreach \c/\r in {0/3,2/3,0/1,2/1,4/2,6/2,4/0,6/0}
      \fill[black!16] (\c,\r) rectangle (\c+1,\r+1);
    \draw[black!45,step=1] (0,0) grid (8,4);
    \draw[thick,densely dashed] (0,1) rectangle (3,4);
    \draw[thick,densely dashed] (4,0) rectangle (7,3);
    \node at (0.5,3.5) {$1$};
    \node at (2.5,1.5) {$2$};
    \node at (4.5,2.5) {$4$};
    \node at (6.5,0.5) {$5$};
    \node at (6.5,2.5) {$3$};
    \draw[->,thick] (0.95,3.5) -- (2.45,3.5);
    \draw[->,thick] (2.05,1.5) -- (0.55,1.5);
    \node[above] at (1.5,4) {\small free};
    \node[above] at (5.5,4) {\small not free};
  \end{scope}
\end{tikzpicture}
\caption{The proof of Lemma~\ref{glob:low} for a tableau of shape
$(N-3,2,1)$, of depth $k=3$. Left: the five marked labels are shaded.
The designated pairs $\{1,2\}$ and $\{4,5\}$ are framed, and $3$ is a
singleton. Right: a configuration of the marked labels on a part of
the board. Every line is light if $n\ge16$. The shaded sites are the
rectangles of the two pairs. The pair $\{1,2\}$ is free, and the
exchange moves its two labels along their rows (arrows). The pair
$\{4,5\}$ is not free, because the singleton lies in its rectangle.}
\label{fig:exchange}
\end{figure}

\step{Step 5: The canceling exchange.}
We exchange the columns of the two labels of a pair. The exchange
must not move a label to a site that carries another marked label.
For alternating line groups the lines involved must also have at
least two unmarked sites, so that the two orbits of Step~4 play no
role. The following definitions, illustrated in
Figure~\ref{fig:exchange}, ensure both.

Call a line \emph{heavy} if it contains more than
$L=\lfloor n^{1/4}\rfloor$ marked labels, and \emph{light} otherwise.
A light line has at least $n-L$ unmarked sites, hence at least two
when $n\ge3$.
The \emph{row set} and the \emph{column set} of a cluster are the sets
of rows and of columns that contain its labels. Its \emph{rectangle}
is the set of sites whose row belongs to the row set and whose column
belongs to the column set. A pair is \emph{free} if the lines of its
row set and of its column set are all light and its rectangle is
disjoint from the rectangle of every other cluster.

Let $\one_{\rm free}$ be the indicator of the event that at least one
pair is free. On this event let $\{i,i'\}$ be the first free pair in a
fixed order of the pairs, and write $\eta(i)=(j,y)$ and
$\eta(i')=(j',y')$. The \emph{exchange} replaces $\eta$ by the map
$\eta^*$ that agrees with $\eta$ on the other marked labels and
satisfies
\[
 \eta^*(i)=(j',y),\qquad \eta^*(i')=(j,y').
\]
Each label keeps its row and takes the column of the other. We claim
that the exchange is an involution of the set of configurations with
a free pair, and that
\[
 F(\eta^*)=F(\eta),\qquad H(\eta^*)=-H(\eta).
\]
The claim implies the cancellation. An involution of a finite set
preserves the uniform law on this set, so
\[
 \E[FH\one_{\rm free}]
 =\E\bigl[F(\eta^*)H(\eta^*)\one_{\rm free}\bigr]
 =-\E[FH\one_{\rm free}],
\]
and therefore
\[
 \E[FH\one_{\rm free}]=0 .
\]

We now prove the claim. The map $\eta^*$ is a configuration. Indeed,
the sites $\eta^*(i)$ and $\eta^*(i')$ are distinct and lie in the
rectangle of the pair. No other marked label lies in this rectangle,
because the pair is free and the site of a marked label belongs to
the rectangle of its own cluster. The exchange does not change the row of any
label, and it interchanges the columns of $i$ and $i'$. The row set,
the column set and the rectangle of every cluster, and the number of
marked labels on every line, are therefore the same for $\eta^*$ as
for $\eta$. So $\eta^*$ has the same free pairs as $\eta$ and the
same first free pair, and the exchange applied to $\eta^*$ gives back
$\eta$. The exchange is thus an involution.

The configurations $\eta^*$ and $\eta$ assign the same row to every
marked label, and they differ only on the rows $y$ and $y'$, which
are light. By Step~4 they lie in the same row orbit, so
$F(\eta^*)=F(\eta)$. Let
$\tau$ be the transposition of $i$ and $i'$, so that
$(\eta\circ\tau)(i)=(j',y')$ and $(\eta\circ\tau)(i')=(j,y)$. The
configurations $\eta^*$ and $\eta\circ\tau$ assign the same column to
every marked label, and they differ only on the columns $j$ and $j'$,
which are light. By Step~4 with rows and columns interchanged they
lie in the same column orbit, and Step~3 gives
$H(\eta^*)=H(\eta\circ\tau)=-H(\eta)$. This argument allows the two
labels to share a row or a column. If they share a column, then
$\eta^*=\eta$, and the claim says that $H(\eta)=0$.

\step{Step 6: The exceptional configurations.}
It remains to control the configurations without a free pair. Let
$\mathcal H_{\rm row}$ be the event that the heavy rows contain more
than $k/16$ marked labels in total, and let $\mathcal H_{\rm col}$ be
the corresponding event for columns. Let $\mathcal N$ be the event
that no pair is free and that neither $\mathcal H_{\rm row}$ nor
$\mathcal H_{\rm col}$ occurs. Every configuration without a free
pair lies in $\mathcal H_{\rm row}$, in $\mathcal H_{\rm col}$ or in
$\mathcal N$. We use two estimates, valid for all sufficiently large
$n$, every row orbit $\mathcal R$ and every column orbit
$\mathcal C$. The first concerns the heavy lines:
\begin{equation}\label{glob:heavy}
 \mathbb P(\mathcal H_{\rm col}\mid\mathcal R)\le n^{-k/200},\qquad
 \mathbb P(\mathcal H_{\rm row}\mid\mathcal C)\le n^{-k/200}.
\end{equation}
The second concerns $\mathcal N$:
\begin{equation}\label{glob:both}
 \mathbb P(\mathcal N\mid\mathcal R)\le n^{-k/40},\qquad
 \mathbb P(\mathcal N\mid\mathcal C)\le n^{-k/40}.
\end{equation}
Both are proved in Appendix~\ref{exc:section}. The first holds
because at most $t/L$ columns can be heavy. For the second, a pair on
light lines that is not free has a label in the same column as a
label of another cluster, and many such coincidences are unlikely.

\step{Step 7: Conclusion.}
For an event $A$, the Cauchy--Schwarz inequality gives
\[
 \E[|FH|\one_A]
 \le \E[F^2\one_A]^{1/2}\E[H^2\one_A]^{1/2}.
\]
By Step~4, $F^2$ is a function of the row orbit and $H^2$ is a
function of the column orbit. We condition on these orbits and recall
that $F$ and $H$ have unit norm. This gives
\[
 \E[F^2\one_A]\le\max_{\mathcal R}\mathbb P(A\mid\mathcal R),\qquad
 \E[H^2\one_A]\le\max_{\mathcal C}\mathbb P(A\mid\mathcal C).
\]
For $A=\mathcal N$, both maxima are at most $n^{-k/40}$ by
\eqref{glob:both}, so $\E[|FH|\one_{\mathcal N}]\le n^{-k/40}$. The
event $\mathcal H_{\rm row}$ is determined by the row orbit, so for
$A=\mathcal H_{\rm row}$ we bound the factor with $H$ by
\eqref{glob:heavy} and the factor with $F$ by $1$. This gives
$\E[|FH|\one_{\mathcal H_{\rm row}}]\le n^{-k/400}$, and in the same
way $\E[|FH|\one_{\mathcal H_{\rm col}}]\le n^{-k/400}$.

By Step~5, $\E[FH]=\E[FH(1-\one_{\rm free})]$. By Step~6, every
configuration without a free pair lies in $\mathcal N$, in
$\mathcal H_{\rm row}$ or in $\mathcal H_{\rm col}$. Therefore
\[
 |\E[FH]|\le n^{-k/40}+2n^{-k/400}\le n^{-k/1000}
\]
for all sufficiently large $n$. This is the bound required in Step~3.
\end{proof}

\begin{remark}[The two smallest shapes]
For $\lambda=(N-1,1)$ the tableau has one designated pair and no other
cluster. For $n\ge16$ this pair is always free, so Step~5 of the proof
gives $\beta_\lambda=0$. For $\lambda=(N-2,2)$ there are two
designated pairs, and neither is free when their rectangles meet. For
this shape the correlation does not vanish. Consider two given labels.
The event that they lie in the same row and the event that they lie
in the same column are disjoint, and each has probability $1/(n+1)$.
Their indicators therefore have correlation $-1/n$. The first
indicator is invariant under the row groups and the second under the
column groups. After centering, both belong to the
$\lambda$-isotypic component, and therefore $\beta_\lambda\ge1/n$.
\end{remark}

\subsection{A fourth trace bound}\label{glob:tracesection}

We bound the remaining correlations through $\operatorname{Tr}(RCRC)$.
For symmetric line groups, we express this trace using common-neighbor
counts in the graph of compatible row and column partitions.
Lemma~\ref{glob:trace} gives the counting bound, and
Lemma~\ref{glob:high} applies it to representations of large depth.

\begin{lemma}\label{glob:trace}
On the full regular representation, for either choice of line groups,
\[
 \operatorname{Tr}(RCRC)\le \exp\!\bigl(O(n\log^2 n)\bigr).
\]
\end{lemma}

\begin{proof}
\step{Step 1: From the trace to compatible partitions.}
First take symmetric line groups. The row data of a board form an
ordered balanced partition of the labels into $n$ sets of size $n$.
Column data have the same form. Two partitions are compatible when
each row block meets each column block exactly once. Each compatible
pair determines precisely one board. Write $M=N!/(n!)^n$ for the
number of partitions, $\Delta=(n!)^n$ for the number of compatible
neighbors of a given partition, and $\kappa=\Delta/M$. Let $K$ have entries
$1/\Delta$ on compatible pairs and zero elsewhere. Each partition is
the row data of exactly $\Delta$ boards, one for each compatible
partition. For a uniform board the row data are therefore uniform on
the $M$ partitions, and given the row data the column data are uniform
on the $\Delta$ compatible partitions. Thus $K$ is the transition
kernel from the row data to the column data, and by symmetry also from
the column data to the row data.

Let $\iota_{\rm row}$ send a function $\varphi$ of the row data to the
function of the board whose value at $g$ is $\varphi$ at the row data
of $g$, and define $\iota_{\rm col}$ in the same way. Both maps are
isometries for the uniform laws. Then
$R=\iota_{\rm row}\iota_{\rm row}^*$ and
$C=\iota_{\rm col}\iota_{\rm col}^*$, because $R$ and $C$ are the
conditional expectations given the row data and the column data, and
$\iota_{\rm row}^*\iota_{\rm col}=K=K^*$ by the previous paragraph.
The cyclicity of the trace gives
\[
 \operatorname{Tr}(RCRC)
 =\operatorname{Tr}\bigl((\iota_{\rm row}^*\iota_{\rm col})
   (\iota_{\rm col}^*\iota_{\rm row})
   (\iota_{\rm row}^*\iota_{\rm col})
   (\iota_{\rm col}^*\iota_{\rm row})\bigr)
 =\operatorname{Tr}(K^4).
\]

The entry $K^2(x,z)$ is $\Delta^{-2}$ times the number of common
neighbors of $x$ and $z$, that is, of partitions compatible with both.
Fix a partition $z$, and let $t_z(x)$ be the probability that a
uniform neighbor of $x$ is compatible with $z$. Then
$K^2(x,z)=t_z(x)/\Delta$, that is, $M K^2(x,z)=t_z(x)/\kappa$. For a
uniform ordered partition $x$,
\begin{equation}\label{glob:traceid}
 \operatorname{Tr}(K^4)=\sum_{x,z'}\bigl(K^2(x,z')\bigr)^2
 =M\sum_x\bigl(K^2(x,z)\bigr)^2=\E_x[(t_z(x)/\kappa)^2].
\end{equation}
The first equality holds because $K^2$ is symmetric. For the second, a
permutation of the labels acts on the partitions, preserves
compatibility and can map any partition to $z$, so the sum over $x$
does not depend on $z'$. The third is the formula for $MK^2(x,z)$.

Write $a_{ij}=|x_i\cap z_j|$ for the overlaps of the blocks of $x$ and
$z$. These multiplicities form a bipartite multigraph whose individual
edges are distinguished by their labels. The blocks of $x$ are the
left vertices, the blocks of $z$ are the right vertices, and each
label is an edge joining the two blocks that contain it. Every vertex
has degree $n$. A common neighbor $w$ of $x$ and $z$ colors each edge
with the index of the block of $w$ that contains its label.
Compatibility with $x$ says that the $n$ edges at each left vertex
have distinct colors, and compatibility with $z$ says the same at each
right vertex. A common neighbor of $x$ and $z$ is therefore precisely a
proper coloring of the edges with $n$ named colors. For example, if
$n=2$, $x=(\{1,2\},\{3,4\})$ and $z=(\{1,3\},\{2,4\})$, the multigraph
is a cycle of length four, and its two proper colorings are the common
neighbors $(\{1,4\},\{2,3\})$ and $(\{2,3\},\{1,4\})$. Let $\Phi(a)$
be the number of proper colorings, which depends only on the overlap
table $a$, so that $t_z(x)=\Phi(a)/\Delta$.

\step{Step 2: Counting common neighbors.}
This step bounds $\Phi(a)$, with the result \eqref{glob:bregman}.
For the overlap table $a$, assign colors one left vertex at a time.
When vertex $i$ is processed,
each edge to $j$ has $s_{ij}=n-\sum_{i'<i}a_{i'j}$ available colors.
The number of assignments of distinct colors to the $n$ edges at $i$
is a zero-one permanent. The rows of its matrix are the $n$ edges at
$i$, the columns are the $n$ colors, and an entry is one when the
color is still available for the edge. Put $g(s)=\log(s!)/s$. The
permanent inequality of Br\`egman \cite{bregman} bounds the permanent
of a zero-one matrix with row sums $r_1,\dots,r_n$ by
$\prod_e(r_e!)^{1/r_e}$. Here $a_{ij}$ rows have sum $s_{ij}$, so the
number of choices is at most
$\exp(\sum_{j:a_{ij}>0}a_{ij}g(s_{ij}))$, independently of the
previous valid assignments.

If every overlap $a_{ij}$ were zero or one, the numbers $s_{ij}$ with
$a_{ij}=1$ would run through $n,n-1,\dots,1$ for each $j$, and the
exponents would add up to $n\sum_{v=1}^ng(v)$. We now bound the error
made when overlaps are larger.
For each $j$ the terms with $a=a_{ij}>0$ and $s=s_{ij}$ correspond
to disjoint integer intervals $[s-a+1,s]$ covering $[1,n]$. Since
$0\le g(v)-g(v-1)\le 1/v$ for $v\ge2$ (use
$\log((v-1)!)\ge (v-1)(\log v-1)$), their individual errors obey
\[
 0\le a g(s)-\sum_{w=s-a+1}^s g(w)
 \le \sum_{v=s-a+2}^s (v-s+a-1)/v\le a .
\]
For intervals of length at most $L_0=\lceil\log n\rceil$ the sum of
these errors at a given $j$ is at most $L_0\log n$, by disjointness.
For the longer intervals we use the bound $a$ of the display. Their
errors add up to at most
$T(a)=\sum_{ij}a_{ij}\one_{\{a_{ij}>L_0\}}$. Since
$\log\Delta=n\log(n!)$, Stirling's estimate
in the form $\sum_{v=1}^n g(v)=\log(n!)-n+O(\log^2 n)$ now yields
\begin{equation}\label{glob:bregman}
 \log(\Phi(a)/\Delta)\le -N+O(n\log^2 n)+T(a).
\end{equation}

\step{Step 3: Averaging over the overlap table.}
To average the error $T(a)$, we use the distribution of the overlap
table. For a uniform $x$, this table has mass
$\kappa/\prod_{ij}a_{ij}!$ on integer tables with all margins $n$.
Indeed, one distributes the $n$ labels at each $z_j$ among the $x_i$,
and divides the product of the resulting multinomial coefficients by
$M$. This is also the law of independent Poisson$(1)$ entries conditional
on all margins being $n$. By Stirling's formula,
$\log\kappa=-N+O(n\log n)$,
so this conditioning has probability $e^{-N}/\kappa=\exp[-O(n\log n)]$.
Before conditioning,
\[
 \E e^{2T}=
 \left(1+e^{-1}\sum_{j>L_0}\frac{e^{2j}-1}{j!}\right)^N=1+o(1).
\]
Indeed, the first term of the tail is $\exp[-(1+o(1))L_0\log L_0]$ and
eventually each successive positive term is at most half its predecessor.
Conditioning on an event of probability $\exp[-O(n\log n)]$ multiplies
the expectation of a nonnegative variable by at most
$\exp[O(n\log n)]$.
Together with \eqref{glob:bregman} this gives
$\E_x t_z(x)^2\le\exp[-2N+O(n\log^2 n)]$. Substitution into
\eqref{glob:traceid} proves the bound for symmetric line groups.

\step{Step 4: Alternating line groups.}
For alternating line groups the matrix entries of either projection,
$R$ or $C$,
are at most $2^n$ times their symmetric counterparts. Indeed, the
average over the alternating group of a line has matrix entries equal
to zero or to twice those of the average over its symmetric group.
The projections $R$ and $C$ are products of $n$ such averages on
disjoint lines.
Expanding the
four-factor trace as a sum of nonnegative products changes its upper
bound by at most $2^{4n}$, which is $\exp(O(n))$.
\end{proof}

The bound of Lemma~\ref{glob:trace} grows with $n$, but on components
of large depth it is small compared with the dimension $d_\lambda$.
The multiplicity $d_\lambda$ in the regular representation then
yields the following dimension-dependent estimate.
\begin{lemma}\label{glob:high}
For all sufficiently large $n$, every $\lambda$ of depth
$k>n^{21/20}$ satisfies $\beta_\lambda\le d_\lambda^{-1/5}$.
\end{lemma}
\begin{proof}
Since $RCR=(RC)(RC)^*$, the trace
$\operatorname{Tr}(RCRC)=\operatorname{Tr}((RCR)^2)$ is the sum of the
fourth powers of the singular values of $RC$. Each singular
value of $RC$ on one copy of $V_\lambda$ occurs at least
$d_\lambda$ times in the regular representation. Hence, by
Lemma~\ref{glob:trace},
\begin{equation}\label{glob:mult}
 d_\lambda\beta_\lambda^4\le \exp[O(n\log^2 n)] .
\end{equation}
It remains to bound $d_\lambda$ from below. We count standard tableaux
in the first case below and use the hook-length formula in the second.
Let $k_0=\lfloor n^{21/20}\rfloor$, and transpose the shape for the
dimension estimate if necessary, so its first row has length
$r=\max(\lambda_1,\lambda'_1)$. Transposition does not change
$d_\lambda$.

\step{Case $r\ge2k_0$.}
The shape has $k>k_0$ boxes below its first row. These boxes form a
Young diagram. Remove corner boxes from it, one at a time, until
$k_0$ boxes remain, and call these the lower boxes. The lower boxes
and the first $2k_0$ boxes of row one form a Young diagram inside
$\lambda$, the subdiagram of Figure~\ref{fig:tableaux}.
Every lower box lies in one of the first $k_0$ columns. Fill the first
$k_0$ boxes of row one with $1,\dots,k_0$. Then choose which $k_0$ of
the next $2k_0$ numbers go to the remaining $k_0$ boxes of row one, in
increasing order. Place the other $k_0$ numbers in the lower boxes,
following a fixed standard ordering of these boxes. Each choice gives
a standard tableau on the subdiagram. Row one increases, and the lower
boxes are filled in a standard order. Every lower box of the second
row holds a number larger than $k_0$, and the box of row one above it
holds a number at most $k_0$. The resulting
$\binom{2k_0}{k_0}$ standard tableaux on the subdiagram all extend to
tableaux on $\lambda$. Hence
$d_\lambda\ge\binom{2k_0}{k_0}\ge2^{k_0}$.

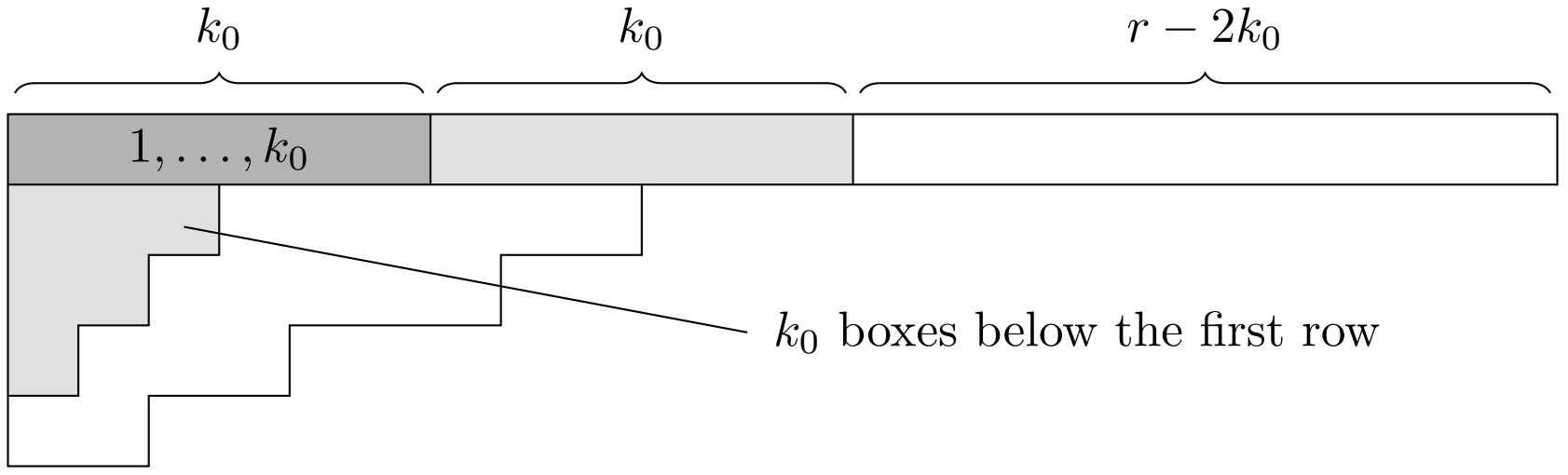
\begin{figure}[ht]
\centering
\begin{tikzpicture}[x=0.5cm,y=0.5cm]
  \fill[black!30] (0,0) rectangle (6,-1);
  \fill[black!12] (6,0) rectangle (12,-1);
  \fill[black!12] (0,-1) -- (3,-1) -- (3,-2) -- (2,-2) -- (2,-3)
    -- (1,-3) -- (1,-4) -- (0,-4) -- cycle;
  \draw (0,0) -- (22,0) -- (22,-1) -- (9,-1) -- (9,-2) -- (7,-2)
    -- (7,-3) -- (4,-3) -- (4,-4) -- (2,-4) -- (2,-5) -- (0,-5) -- cycle;
  \draw (0,-1) -- (9,-1);
  \draw (6,0) -- (6,-1);
  \draw (12,0) -- (12,-1);
  \draw (3,-1) -- (3,-2) -- (2,-2) -- (2,-3) -- (1,-3) -- (1,-4) -- (0,-4);
  \node at (3,-0.5) {\small $1,\dots,k_0$};
  \draw[decorate,decoration={brace,amplitude=4pt}]
    (0.1,0.3) -- (5.9,0.3) node[midway,above=5pt] {\small $k_0$};
  \draw[decorate,decoration={brace,amplitude=4pt}]
    (6.1,0.3) -- (11.9,0.3) node[midway,above=5pt] {\small $k_0$};
  \draw[decorate,decoration={brace,amplitude=4pt}]
    (12.1,0.3) -- (21.9,0.3) node[midway,above=5pt] {\small $r-2k_0$};
  \node[anchor=west] at (10.6,-3.1) {\small $k_0$ boxes below the first row};
  \draw (10.5,-3.1) -- (2.5,-1.6);
\end{tikzpicture}
\caption{The subdiagram used in the proof of Lemma~\ref{glob:high}
when $r\ge2k_0$, shaded inside the shape $\lambda$. It consists of the
first $2k_0$ boxes of the first row and of $k_0$ boxes below it, which
lie in the first $k_0$ columns. The numbers $1,\dots,k_0$ fill the
dark boxes. The next $2k_0$ numbers are divided between the other
shaded boxes of the first row and the lower shaded boxes.}
\label{fig:tableaux}
\end{figure}

\step{Case $r<2k_0$.}
Every hook has length at most
$4k_0$, and the hook-length formula gives
$d_\lambda\ge N!/(4k_0)^N\ge(N/(4ek_0))^N$.

Thus $\log d_\lambda\ge c k_0$ in both cases, uniformly in the stated
range. In \eqref{glob:mult},
$n\log^2 n=o(\log d_\lambda)$, so that
$\beta_\lambda\le d_\lambda^{-1/4+o(1)}$, which proves the claim.
\end{proof}

\subsection{Interpolation and entropy duality}\label{glob:interpolation}

We use the two spectral bounds to prove the centered $L^p$ estimate,
then restore the means and apply entropy duality.

\begin{proof}[Proof of Theorem~\ref{thm:global}]
\step{Step 1: Interpolation on each component.}
Until Step~3 we work on $S_N$ with its uniform probability, as in
Subsection~\ref{glob:strategy}. Let $E_\lambda$ be the orthogonal
projection of $L^2(S_N)$ onto the isotypic component for $\lambda$.
Its rank is $d_\lambda^2$. Write $e_\lambda(g,g')$ for the kernel of
$E_\lambda$ with respect to the uniform probability, so that
$E_\lambda f(g)=\E[e_\lambda(g,g')f(g')]$ for a uniform $g'$. The mean
of the diagonal values $e_\lambda(g,g)$ is the trace of $E_\lambda$,
which equals its rank, and these values do not depend on $g$ by
transitivity. Hence $e_\lambda(g,g)=d_\lambda^2$ for every $g$. Since
$E_\lambda^2=E_\lambda=E_\lambda^*$, the squared $L^2$ norm of every
kernel row $e_\lambda(g,\cdot)$ is this same diagonal value. By the
Cauchy--Schwarz inequality, $|E_\lambda f(g)|\le d_\lambda\|f\|_2$ for
every $g$, with equality for $f=e_\lambda(g,\cdot)$. Thus
$\|E_\lambda\|_{2\to\infty}=d_\lambda$, and by duality
$\|E_\lambda\|_{1\to2}=d_\lambda$. Interpolation between this bound
and the $L^2$ contraction $\|E_\lambda\|_{2\to2}\le1$ gives
$\|E_\lambda\|_{p\to2}\le d_\lambda^\theta$, where $\theta=2/p-1$.

Write $P_0$ for the projection onto the constants for symmetric line
groups and onto the
functions of global parity for alternating line groups. Each $E_\lambda$
commutes with $R$ and $C$. For invariant functions $F=RF$ and $H=CH$, with
norms on the full group, the orthogonality of the isotypic components
gives
$\langle F,H\rangle=\sum_\lambda\langle E_\lambda F,E_\lambda H\rangle$.
For a shape of depth $k\ge1$, the functions $E_\lambda F=RE_\lambda F$
and $E_\lambda H=CE_\lambda H$ lie in the ranges of $R$ and $C$ on the
$\lambda$-component. The definition of $\beta_\lambda$ and the bound
on $\|E_\lambda\|_{p\to2}$ therefore give
\[
 |\langle E_\lambda F,E_\lambda H\rangle|
 \le\beta_\lambda\|E_\lambda F\|_2\|E_\lambda H\|_2
 \le d_\lambda^{2\theta}\beta_\lambda\|F\|_p\|H\|_p .
\]
The two shapes of depth zero contribute $\langle F,P_0H\rangle$ in
total, as we check after the next display. Summing over the shapes
gives
\begin{equation}\label{glob:interp}
 |\langle F,H\rangle-\langle F,P_0H\rangle|
 \le \|F\|_p\|H\|_p \sum_{\lambda:\ k\ge1}
                 d_\lambda^{2\theta}\beta_\lambda .
\end{equation}
The two shapes with $k=0$ are the trivial and sign representations.
In the symmetric case the sign component of each invariant function
vanishes. In the alternating case $P_0$ includes both shapes.

\step{Step 2: Summing the components.}
Every shape of depth $k$ satisfies $d_\lambda\le N^k$. Indeed, after
transposing if needed, the Young branching rule shows that the
restriction of $V_\lambda$ to the stabilizer $S_{N-k}$ contains the
trivial representation. This places $V_\lambda$ in the permutation
module of ordered $k$-tuples, whose dimension is at most $N^k$. At
each depth there are at most $2^{k+1}$ shapes. Indeed, after
transposing if needed, a shape of depth $k$ is determined by the
partition of $k$ that remains when its first row is deleted, and $k$
has at most $2^k$ partitions. The low-depth
contribution in \eqref{glob:interp} is therefore at most
\[
 \sum_{k\ge1}2^{k+1}n^{(4\theta-1/1000)k}=o(1)
 \qquad(4\theta<1/1000)
\]
by Lemma~\ref{glob:low}. There are $\exp(O(\sqrt N))$ partitions of
$N$ in total. This follows by bounding the logarithm of the partition
generating product at $e^{-z}$ by
$\sum_{j\ge1}[j(e^{zj}-1)]^{-1}\le\pi^2/(6z)$ and taking $z=N^{-1/2}$.
By Lemma~\ref{glob:high} and its dimension estimate, the high-depth
sum is at most $\exp[O(n)-(1/5-2\theta)c n^{21/20}]=o(1)$ when
$2\theta<1/5$. For example, $\theta=1/8000$ satisfies both conditions.

\step{Step 3: Passage to a single class and to all board sizes.}
For symmetric line groups there is one class, and $\pi$ is the uniform
probability on $S_N$.
In the alternating case extend class-centered $F$ and $H$ by zero on
the other parity class. Their extensions are still invariant and
orthogonal to the constants and to the sign character, so the term
with $P_0$ in
\eqref{glob:interp} vanishes for them. Each class has probability
$1/2$ in $S_N$. The inner product of the extensions is therefore half
of $\E_\pi[FH]$, and their $p$-norms on $S_N$ are $2^{-1/p}$ times the
$p$-norms on the class. The conversion from \eqref{glob:interp}
to probability norms on one class thus introduces the factor
$2^{1-2/p}\le1$.

For $p<2$ let $c_n(p)$ be the smallest constant such that, on every
class of the board of side $n$,
\[
 |\E_\pi[FH]|\le c_n(p)\|F\|_p\|H\|_p
 \quad \text{if }\E_\pi F=\E_\pi H=0,\quad F=RF,\ H=CH.
\]
Let $p_0<2$ be the exponent with $\theta=1/8000$. Steps~1 and~2 and
the preceding paragraph show that $c_n(p_0)=o(1)$. Fix $n_0$ such that
$c_n(p_0)\le p_0-1$ for all $n\ge n_0$.

Consider one of the finitely many remaining sizes $n<n_0$. The
centered $L^2$ ranges of $R$ and $C$ have no nonzero common vector.
Indeed, a function in both ranges is invariant under every line group,
and it is therefore constant on the class, which is an orbit of the
group that the line groups generate. The maximal
correlation $c_n(2)$ of the two ranges is strictly less than one by compactness of the
unit spheres. On a finite probability space the $p$-norm is continuous
in $p$, uniformly on bounded sets of functions, so $c_n(p)\to c_n(2)$
as $p\to2$. Since $p-1\to1$, there is $p_n<2$ with
$c_n(p_n)\le p_n-1$.

On a probability space the $p$-norm does not decrease when $p$
increases. Hence $c_n(p)$ does not increase with $p$, whereas $p-1$
increases. The inequality $c_n(p)\le p-1$, once true for some $p<2$,
therefore holds for every larger $p<2$. Let $p$ be the largest of
$p_0$ and of the $p_n$ with $n<n_0$. Then $p<2$, and
$c_n(p)\le p-1$ for all $n$.

\step{Step 4: Restoring the means.}
The factor $p-1$ matches the following convexity inequality and gives
coefficient one in \eqref{glob:bilinear}, as required for entropy duality.
For a real $g$ on a finite probability space,
\begin{equation}\label{glob:convex}
 \|g\|_p^2\ge(\E_\pi g)^2+(p-1)\|g-\E_\pi g\|_p^2,\qquad 1<p<2.
\end{equation}

To verify \eqref{glob:convex}, put $v=g-\E_\pi g$ and
$g_t=\E_\pi g+t v$.
Away from zeros, differentiating $\|g_t\|_p^2$ twice shows that its
second derivative is at least
$2(p-1)\|g_t\|_p^{2-p}\E_\pi[|g_t|^{p-2}v^2]\ge2(p-1)\|v\|_p^2$,
the last step being H\"older's inequality. Integrating from $t=0$,
where the first derivative vanishes, gives \eqref{glob:convex}.
Zero coordinates are handled by continuity, or by integrating between the
finitely many zeros, across which the first derivative is continuous.

Now let $F=RF$ and $H=CH$ be nonnegative, and put $F_0=F-\E_\pi F$ and
$H_0=H-\E_\pi H$. These functions are centered and satisfy $F_0=RF_0$
and $H_0=CH_0$, so the estimate of Step~3 applies to them. Hence
\begin{align*}
 \E_\pi[FH]&=\E_\pi F\,\E_\pi H+\E_\pi[F_0H_0]\\
 &\le\E_\pi F\,\E_\pi H+(p-1)\|F_0\|_p\|H_0\|_p\\
 &\le\bigl((\E_\pi F)^2+(p-1)\|F_0\|_p^2\bigr)^{1/2}
     \bigl((\E_\pi H)^2+(p-1)\|H_0\|_p^2\bigr)^{1/2}
 \le\|F\|_p\|H\|_p .
\end{align*}
The second line is the estimate of Step~3. The third line is the
Cauchy--Schwarz inequality in $\mathbb R^2$, followed by
\eqref{glob:convex} for $F$ and for $H$. This is
\eqref{glob:bilinear}.

\step{Step 5: Passing to entropy.}
By homogeneity we may take $\E_\pi h=1$. Apply \eqref{glob:bilinear} to
$(Rh)^{1/p}$ and $(Ch)^{1/p}$, whose $p$-norms are both one. For
$\phi=p^{-1}(\log Rh+\log Ch)$ this says $\E_\pi[e^\phi]\le1$. The
variational formula for relative entropy and the conditional-expectation
identities for $R$ and $C$ yield
\[
 \Ent_\pi(Rh)+\Ent_\pi(Ch)=p\,\E_\pi[h\phi]\le p\Ent_\pi(h).
\]
This is the first inequality in \eqref{glob:factorization}. Zero
values can be handled by applying the argument to
$(h+\varepsilon)/(1+\varepsilon)$ and taking the limit. For the second, recall
that a row orbit is an orbit of the product of the row groups. The
uniform law on a row orbit is the product of the uniform laws on the
orbits of the individual row groups. Conditional entropy tensorization on
that orbit, followed by integration over orbits, gives
\[
 \Ent_\pi(h)-\Ent_\pi(Rh)
 \le \sum_{\ell\text{ a row}}[\Ent_\pi(h)-\Ent_\pi(B_\ell h)] .
\]
Adding the column version proves \eqref{glob:factorization}.
\end{proof}

\section{Comparison with unit rotations}\label{phys:section}

This section proves the second estimate in \eqref{eq:bridge}. We
measure the cost of a physical permutation by the following energy.
For a permutation $s$ of the sites preserving the communicating class, write
$U_sf(x)=f(s^{-1}x)$ and $e_s(f)=\|(I-U_s)f\|_{2,\pi}^2$.
Then $e_s=e_{s^{-1}}$. Each $U_s$ is an isometry of $L^2(\pi)$ and
$U_{st}=U_sU_t$, so the triangle inequality gives
\begin{equation}\label{phys:triangle}
  \sqrt{e_{st}(f)}\le\sqrt{e_s(f)}+\sqrt{e_t(f)} .
\end{equation}
For a row $y$, let $R_y$ be the permutation of the sites that maps
$(j,y)$ to $(j+1,y)$ for every $j$ and fixes the sites outside the
row. For a column $j$, let $C_j$ map $(j,y)$ to $(j,y+1)$ for every
$y$ and fix the sites outside the column. We call $R_y$ and $C_j$ the
positive unit rotations of these lines. The normalization of the chain
gives
\begin{equation}\label{phys:normalization}
  8n\cE_P(f,f)=\sum_{y}e_{R_y}(f)+\sum_{j}e_{C_j}(f),
\end{equation}
where $y$ runs over the rows and $j$ over the columns.

\begin{theorem}\label{thm:physical}
There is an absolute $C$ such that on each communicating class, for all
$n\ge2$ and all $h\ge0$,
\[
 D_{\rm line}(h)\le C n^2\cE_P(\sqrt h,\sqrt h).
\]
\end{theorem}

\begin{remark}[The factor $n^2$]
The factor $n^2$ has the correct order. For $n\geq3$ and
$0<\varepsilon<1$ let $h=(1+\varepsilon\cos(2\pi u/n))^2$, where $u$
is the horizontal coordinate of a fixed card. Only the row averages
change $h$. A row contains the card with probability $1/n$, and its
average removes this fraction of the entropy of $h$. Hence
$D_{\rm line}(h)=\Ent_\pi(h)/(2n)$. By the proof of
Lemma~\ref{lem:observable}(ii), and since $\cos^2(2\pi u/n)$ has mean
$1/2$, we have
$\cE_P(\sqrt h,\sqrt h)=\varepsilon^2(1-\cos(2\pi/n))/(8n)$. As in
Step~4 of Subsection~\ref{subsec:lower},
$\Ent_\pi(h)=\varepsilon^2+O(\varepsilon^3)$ as $\varepsilon\to0$ with
$n$ fixed. The ratio $D_{\rm line}(h)/(n^2\cE_P(\sqrt h,\sqrt h))$
therefore tends to $4/(n^2(1-\cos(2\pi/n)))\geq2/\pi^2$.
\end{remark}

\paragraph{Strategy of the comparison.}
The line average rearranges all the labels of a line, whereas an
update of the shuffle rotates the whole line. The rotations of one
line produce only its $n$ cyclic shifts, so the labels of a row can be
rearranged only with the help of the columns. Displacing the row by
$d$ sites, rotating a column, and undoing both moves cycles three
sites, two of them on the row (Step~1 in the proof of
Lemma~\ref{phys:routing}). By \eqref{phys:triangle} the displacement
costs at most $d^2$ times the energy of the unit rotation. Since
labels have to travel distances of order $n$, this accounts for the
factor $n^2$.

We bound the entropy
loss under the line average by the energy of the rotations,
recursively bisecting the line.
At each division of an interval of
length $m$ into two halves, its children, the entropy chain rule
separates the child label sets
from their internal orders. We estimate the label-set entropy at the
current scale and apply the same decomposition to each child. A parity
bit at the cut is handled by a bounded-size move.

Write $f=\sqrt h$ and fix a target row $y$. On an interval $I$ of length $m$, we use two
overlapping layers of rearrangements on intervals of length comparable
to a scale $b$.
Lemma~\ref{phys:split} controls the label-set entropy by
$C((m/b)^2+m)$ times the sum of the two layer projection losses.
Lemma~\ref{phys:routing} bounds each layer's projection loss by
$C(b^2e_{R_y}(f)+\sum_{j\in I}e_{C_j}(f))$, where $R_y$ rotates the target
row and only columns meeting the current interval $I$ are used.
Taking $b=\lfloor\sqrt m\rfloor$ makes $(m/b)^2+m$ of order $m$.
A larger $b$ would increase the coefficient of the row rotation in
the product of the two bounds, and a smaller $b$ that of the columns.
Thus the contribution of a node has the form
$C(m^2e_{R_y}(f)+m\sum_{j\in I}e_{C_j}(f))$; this is
\eqref{phys:node} below.

The squared interval lengths have sum $O(n^2)$ over the binary tree.
For a fixed column, only the intervals on one branch contribute,
because the column sum is restricted to $I$. Their lengths have sum
$O(n)$. The resulting single-row bound is a constant times
$n^2e_{R_y}(f)+n\sum_j e_{C_j}(f)$, without a loss for the number of scales.
Averaging this bound over rows and columns gives the theorem.

Subsection~\ref{phys:gapsection} constructs a local update by
disjoint three-cycles with a uniform spectral gap.
Subsection~\ref{phys:realization} realizes it by rotations.
Subsection~\ref{phys:splitsection} proves the label-set
entropy estimate, and Subsection~\ref{phys:proofsection} assembles
the proof of Theorem~\ref{thm:physical}.

\subsection{A spectral gap for disjoint triples}\label{phys:gapsection}

We compare complete rearrangements on short intervals with products
of disjoint three-cycles. The cycles have common spacings,
which allow them to share the same row rotations in the construction
of Lemma~\ref{phys:routing}. We prove that this update has a uniform
spectral gap even with the spacings shared across all intervals.

Fix a family $\mathcal I$ of disjoint consecutive intervals in a
single row. Read
each interval in the direction of the positive rotation of the row
and, when an
interval is the whole row, choose a starting site. Suppose their sizes $m_i$
lie in $[b/2,2b]$ for an integer $b\ge8$. Write $B_I^+$ for uniform
averaging over the alternating group on $I$ and
$B_{\mathcal I}^+=\prod_{I\in\mathcal I}B_I^+$. The random update is
defined in four steps.
\begin{enumerate}
\item Choose two independent
uniform integers $d,d'$ from $\{1,\ldots,2b\}$.
\item In each interval form all
triples of sites $(j-d,j,j+d')$ in its linear order. These are the
candidates.
\item Conditional on $d$ and $d'$,
give each candidate an independent uniform $[0,1]$ priority. Retain a
candidate if its priority is the smallest among all candidates meeting
it, including itself. The retained triples form a disjoint matching
$\mathcal M$.
\item On each retained triple independently, apply a uniform
cyclic permutation, including the identity.
\end{enumerate}
Let $K_{\mathcal I}$ denote
the average operator of this update. For self-adjoint operators $A$
and $B$ we write $A\succeq B$ when $A-B$ is positive semidefinite.

\begin{lemma}\label{phys:gap}
In every unitary representation of the product of the interval
alternating groups,
\[
 I-K_{\mathcal I}\succeq c(I-B_{\mathcal I}^+)
\]
for an absolute $c>0$.
\end{lemma}
\begin{proof}
\step{Step 1: One interval; the matching.}
Consider one interval of length $m$. Given an unordered triple $T$, its
successive spacings specify a pair $d,d'$, which is chosen with
probability $(2b)^{-2}$. Conditional on those spacings, at most
nine candidates meet $T$, since any site belongs to at most three
candidates. The priorities are exchangeable with no ties. Hence
\begin{equation}\label{phys:include}
 \mathbb P(T\in\mathcal M)\ge(36b^2)^{-1}.
\end{equation}
For either cyclic orientation $\tau$ of $T$ put
$p_T=(I+U_\tau+U_{\tau^{-1}})/3$. Conditional on $\mathcal M$, the
update operator is the orthogonal projection $\prod_{T\in\mathcal M}p_T$.

\step{Step 2: Simultaneous and individual projections.}
We compare the simultaneous projection loss with the sum of the
individual triple losses on each Specht module of $S_m$. Let $\lambda$
index this module and put
$r=\max(\lambda_1,\lambda'_1)$ and $k=m-r$, the depth of $\lambda$ as
in Section~\ref{glob:section}. The modules with $k=0$ are
trivial upon restriction to $A_m$. For $k>0$ we may transpose the
diagram so that $\lambda_1=r$, because conjugate diagrams differ by a
sign twist. The branching rule places this module in the permutation
representation on ordered $k$-tuples of distinct sites because removing the
$k$ boxes below the first row exhibits a trivial $S_{m-k}$ constituent
upon restriction. For fixed $\mathcal M$, the cyclic groups of the
triples in $\mathcal M$ commute, because the triples are disjoint.
Consider the orbit of any
$k$-tuple under the product of these groups. Every triple missing
the tuple acts trivially throughout its orbit, and at most $k$ triples
of $\mathcal M$ meet the tuple. Decompose the functions on this orbit
into joint eigenfunctions of the commuting projections $p_T$ with
$T\in\mathcal M$, whose eigenvalues are zero or one. On a joint
eigenfunction at most $k$ of the $p_T$ can therefore vanish. If none
vanishes, both sides of \eqref{phys:active} below are zero on it.
Otherwise the left-hand side equals one and the right-hand side is at
most one.
Summing over the orbits and restricting to the chosen
Specht module gives
\begin{equation}\label{phys:active}
 I-\prod_{T\in\mathcal M}p_T
 \succeq k^{-1}\sum_{T\in\mathcal M}(I-p_T).
\end{equation}

\step{Step 3: The sum over all triples.}
After averaging \eqref{phys:active}, the inclusion bound
\eqref{phys:include} reduces the problem to a sum over all unordered
triples. We show that its eigenvalue is at least $ckm^2$, which cancels
the factor $k^{-1}$ in \eqref{phys:active}. The sum of $I-p_T$ is
central, and its eigenvalue on $\lambda$ is
\begin{equation}\label{phys:content}
 \frac13\left(\sum_{j=0}^{m-1}j^2-
                    \sum_{z\in\lambda}{\rm cont}(z)^2\right),
\end{equation}
where the content of a box is its column minus its row. Indeed,
$J_i=\sum_{j<i}(j\,i)$ satisfies $\sum_iJ_i^2=\binom m2 I+
\sum_{\tau}\tau$, with the last sum over all oriented three-cycles.
In the tableau basis $J_i$ acts by the content of the box containing $i$
\cite{vershikokounkov}. Subtracting the eigenvalue for $\lambda$ from the
one-row value gives \eqref{phys:content} after division by three.

For $k\le m/3$, the missing tail of the one-row sum contributes at least
$k(m-k)^2$, and the $k$ boxes below the first row each have absolute
content at most $k$. Thus \eqref{phys:content} is at least $c k m^2$.
For $k>m/3$, the $i$th box in any standard tableau has absolute content
at most $\min(i-1,r-1)$. Since $r<2m/3$, for a fixed positive fraction
of indices near $m$ the deficit from $(i-1)^2$ is at least $c m^2$.
The same conclusion holds for the finitely many small $m$ after changing
$c$, since three-cycles generate $A_m$. Average \eqref{phys:active} and
use \eqref{phys:include}, \eqref{phys:content} and $m\ge b/2$ to obtain
a uniform spectral gap on every nontrivial factor: on such a module
the averaged update operator $K$ of the interval satisfies
$I-K\succeq(36b^2k)^{-1}ckm^2\ge c/144$.

\step{Step 4: Products of intervals.}
The shared spacings $d$ and $d'$ make the interval updates dependent.
We can still
apply the one-interval bound by retaining the loss from a single
nontrivial factor. Decompose the regular representation of
$\prod_i S_{m_i}$ into tensor products of Specht modules.
On any component with a factor
nontrivial on $A_{m_i}$, the complement of the full matching projection
dominates the complement of its projection on this factor tensored with
the identity on the other factors. Taking expectations uses only the
one-interval marginal, so shared spacings cause no loss. On components
trivial on each alternating group both sides vanish. Restriction of the
regular action to $\prod_i A_{m_i}$ contains its regular representation,
which contains every irreducible representation. This proves the
lemma.
\end{proof}

\subsection{Realization by row and column rotations}\label{phys:realization}

We implement the three-cycle layers of Lemma~\ref{phys:gap} using row
displacements of order $b$ and column displacements of one unit.
The estimate below keeps the column energies restricted to the
intervals being rearranged, as required for the recursive sum.

\begin{lemma}\label{phys:routing}
Let $R_y$ rotate the target row $y$ by one unit and $C_j$ rotate column $j$.
For intervals as in Lemma~\ref{phys:gap},
\begin{equation}\label{phys:routing-est}
 \|(I-B_{\mathcal I}^+)f\|_{2,\pi}^2
 \le C\left(b^2 e_{R_y}(f)+
        \sum_{j\in\bigcup_{I\in\mathcal I} I}e_{C_j}(f)\right).
\end{equation}
In the sum we identify sites of the target row with their column indices.
\end{lemma}
\begin{proof}
\step{Step 1: Double swaps.}
Recall that products of permutations are composed from right to
left. Lists of physical moves below are in
chronological order. We first exchange pairs on the target row while
allowing an auxiliary swap off the row. Common spacings let us
implement all these exchanges with the same row displacements. We
remove the auxiliary swaps after assembling the desired three-cycles.

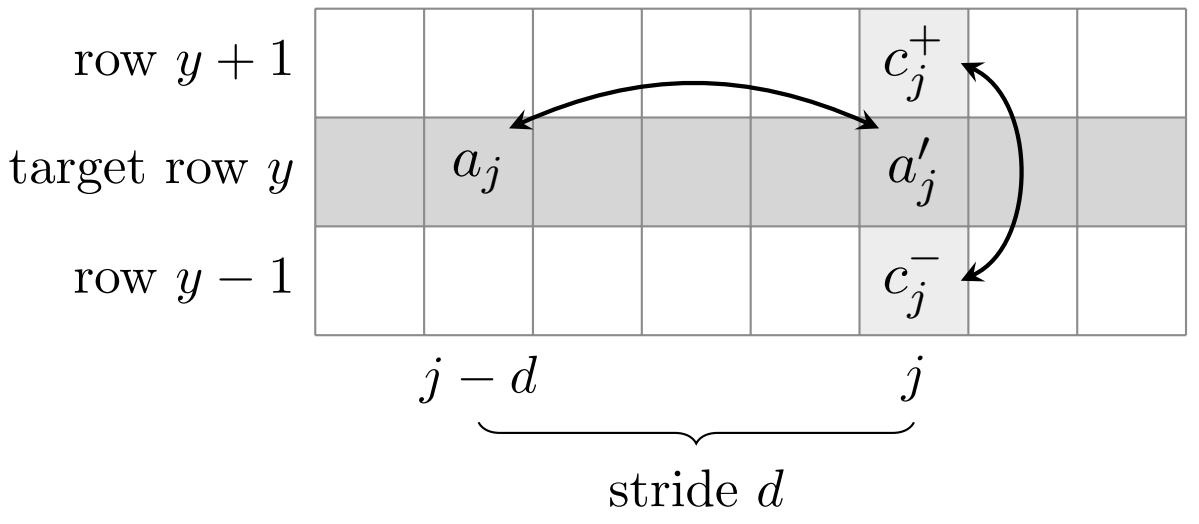
\begin{figure}[ht]
\centering
\begin{tikzpicture}[x=0.72cm,y=0.72cm,>=stealth]
  \fill[black!7]  (5,0) rectangle (6,3);
  \fill[black!16] (0,1) rectangle (8,2);
  \draw[black!45,step=1] (0,0) grid (8,3);
  \node at (1.5,1.5) {$a_j$};
  \node at (5.5,1.5) {$a'_j$};
  \node at (5.5,2.5) {$c_j^+$};
  \node at (5.5,0.5) {$c_j^-$};
  \node[left]  at (0,2.5) {\small row $y+1$};
  \node[left]  at (0,1.5) {\small target row $y$};
  \node[left]  at (0,0.5) {\small row $y-1$};
  \node[below] at (1.5,0) {\small $j-d$};
  \node[below] at (5.5,0) {\small $j$};
  \draw[<->,thick] (1.78,1.9) to[bend left=24] (5.18,1.9);
  \draw[<->,thick] (5.93,2.5) to[bend left=65] (5.93,0.5);
  \draw[decorate,decoration={brace,mirror,amplitude=4pt}]
    (1.5,-0.8) -- (5.5,-0.8) node[midway,below=5pt] {\small stride $d$};
\end{tikzpicture}
\caption{One quadruple in the proof of Lemma~\ref{phys:routing}. The
sites $a_j=(j-d,y)$ and $a'_j=(j,y)$ lie on the target row, and
$c_j^\pm=(j,y\pm1)$ are the neighbors of $a'_j$ in column~$j$. The
double swap $D_J$ exchanges $a_j$ with $a'_j$ and $c_j^+$ with $c_j^-$
(arrows). It is built from the row displacements $R_y^{\pm d}$ and from
unit rotations of column~$j$.}
\label{fig:doubleswap}
\end{figure}

Let the row be $y$ and consider disjoint pairs
$a_j=(j-d,y)$ and $a'_j=(j,y)$ for $j\in J$, all with the same stride $d$.
Put $c_j^\pm=(j,y\pm1)$; see Figure~\ref{fig:doubleswap}. We first construct the deterministic permutation
$D_J=\prod_{j\in J}(a_j\,a'_j)(c_j^+\,c_j^-)$ and show that
\begin{equation}\label{phys:double}
 \sqrt{e_{D_J}(f)}
 \le C\left(d\sqrt{e_{R_y}(f)}+
                 \sqrt{\sum_{j\in J}e_{C_j}(f)}\right).
\end{equation}
Use independent uniform signs~$\varepsilon_j$ and
let~$Z=\prod_{j\in J}C_j^{\varepsilon_j}$.
Applying~$R_y^d,Z,R_y^{-d},Z^{-1}$ in that order gives
\begin{equation*}
 \omega_+=Z^{-1}R_y^{-d}ZR_y^d=\prod_j(a_j\,a'_j\,c_j^{-\varepsilon_j})\,.
\end{equation*}
For a single column~$j$, the four moves send~$a_j$ to~$a'_j$, $a'_j$
to~$c_j^{-\varepsilon_j}$, and~$c_j^{-\varepsilon_j}$ to~$a_j$,
fixing every other site. For each selected index the conjugated
column has the same off-row sites as the original column and has its
row site replaced by~$a_j$. Thus the selected disjoint endpoints
prevent any interaction between columns. Reversing the signs
gives~$\omega_-=\prod_j(a_j\,a'_j\,c_j^{\varepsilon_j})$.

The two cycles share the row edge and use opposite off-row sites. On
one quadruple, suppress~$j$ and take~$\varepsilon_j=1$. The images of
the four positions under the successive factors are
\begin{equation*}
 \begin{array}{c|cccc}
  x & a & a' & c^- & c^+ \\
  \hline
  \omega_+(x) & a' & c^- & a & c^+ \\
  \omega_-^2\omega_+(x) & a & c^- & c^+ & a' \\
  \omega_+\omega_-^2\omega_+(x) & a' & a & c^+ & c^-
 \end{array}\,.
\end{equation*}
The final row is the double swap~$(a\,a')(c^+\,c^-)$.
For~$\varepsilon_j=-1$, interchange~$c^+$ and~$c^-$ in the
calculation. Hence~$\omega_+\omega_-^2\omega_+=D_J$ for every sign choice.

\step{Step 2: Energy of a double swap.}
The permutation $D_J$ is independent of the signs. We average over
them to bound the cost of simultaneous column rotations by the sum
of the individual column energies.
The column rotations commute and admit a joint character
decomposition. If their
multipliers are $e^{{\rm i}\theta_j}$, the multiplier of the average of
$Z$ is $\prod_{j\in J}\cos\theta_j$. Since
$1-\prod_j t_j\le\sum_j(1-t_j)$ for $-1\le t_j\le1$,
\[
 \E_{\varepsilon} e_{Z^{\pm1}}(f)\le\sum_{j\in J} e_{C_j}(f).
\]
The row displacement costs at most $d\sqrt{e_{R_y}(f)}$ in $L^2$,
by \eqref{phys:triangle}. Telescope
the word for $D_J$ in the joint $L^2$ space of the board and the signs.
Conditional on
the entire sign vector, each prefix preserves the uniform board measure.
The preceding estimate therefore applies to every occurrence of $Z^{\pm1}$
and to the reversed sign vector. The word has a bounded number of such
occurrences, proving \eqref{phys:double}.

\step{Step 3: Cancellation off the row.}
Fix a matching of triples and a prescribed set of positive
orientations. For a selected
triple~$(x_1,x_2,x_3)=((j-d,y),(j,y),(j+d',y))$, take the positive
orientation to be~$(x_1\,x_2\,x_3)$. Apply the double-swap layer for
the edges~$(x_2,x_3)$ first and that for the edges~$(x_1,x_2)$
second. The strides are~$d'$ and~$d$, respectively, and the row
permutation on each triple
is~$(x_1\,x_2)(x_2\,x_3)=(x_1\,x_2\,x_3)$. Thus the full row
permutation~$\Gamma$ is a product of disjoint positive three-cycles.

The off-row permutation~$\Theta$ is a product of commuting
transpositions~$(c_j^+\,c_j^-)$, at the columns of the second and
third sites of each triple. Since~$\Gamma^3=\Theta^2=I$ and their
supports are disjoint, four repetitions
give~$(\Gamma\Theta)^4=\Gamma^4\Theta^4=\Gamma$. Use the inverse
construction on the triples assigned the opposite orientation and
omit inactive triples. This realizes every activity and orientation
pattern with a bounded number of double-swap layers, with strides at
most~$2b$ and all column indices in the stated union.

\step{Step 4: Conclusion.}
Let $\sigma$ be the permutation applied by one realization of the
random update of Subsection~\ref{phys:gapsection}, that is, by a
matching and a cyclic permutation of each retained triple. By Step~3,
$\sigma$ is a product of a bounded number of double-swap layers.
Telescoping with \eqref{phys:triangle} and applying
\eqref{phys:double} to each layer bounds the squared increment
$e_\sigma(f)$ by the
right-hand side of \eqref{phys:routing-est}. Since $U_\sigma$ is an
isometry and $K_{\mathcal I}$ is the mean of $U_\sigma$, the
expectation of $e_\sigma(f)$ over the matching
and orientations is $2\langle f,(I-K_{\mathcal I})f\rangle_\pi$.
Lemma~\ref{phys:gap}, applied to the representation of the interval
alternating groups on $L^2(\pi)$, gives
$\langle f,(I-K_{\mathcal I})f\rangle_\pi
\ge c\|(I-B_{\mathcal I}^+)f\|_{2,\pi}^2$, because $B_{\mathcal I}^+$
is an orthogonal projection. This proves \eqref{phys:routing-est}.
\end{proof}

\subsection{The entropy of a balanced division}\label{phys:splitsection}

We now estimate the entropy of the label set in one half of a line
interval. Two layers of local averaging separate the centered input
into a smooth part with controlled entropy and a remainder with small
$L^2$ norm. The remainder is bounded using the number of half-label
sets, whose logarithm is of order $m$. The internal orders within
the halves are left for the recursion.

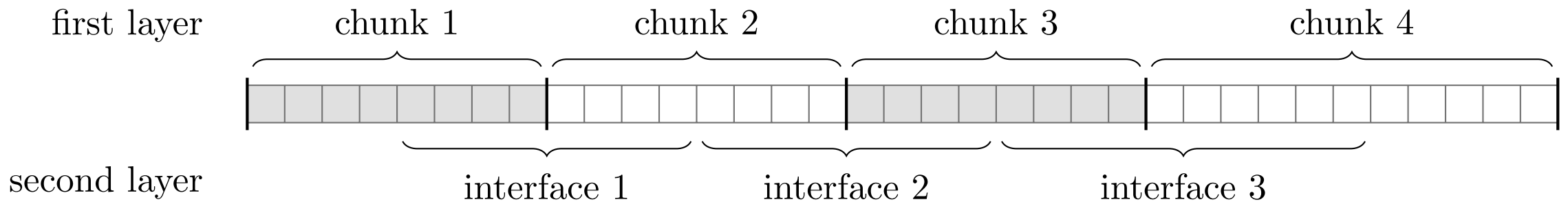
\begin{figure}[ht]
\centering
\begin{tikzpicture}[x=0.36cm,y=0.36cm]
  \fill[black!12] (0,0) rectangle (8,1);
  \fill[black!12] (16,0) rectangle (24,1);
  \draw[black!55,step=1] (0,0) grid (35,1);
  \foreach \x in {0,8,16,24,35} \draw[thick] (\x,-0.2) -- (\x,1.2);
  \foreach \a/\b/\k in {0/8/1,8/16/2,16/24/3,24/35/4}
    \draw[decorate,decoration={brace,amplitude=4pt}]
      (\a+0.15,1.5) -- (\b-0.15,1.5)
      node[midway,above=5pt] {\small chunk \k};
  \foreach \a/\b/\k in {4/12/1,12/20/2,20/30/3}
    \draw[decorate,decoration={brace,mirror,amplitude=4pt}]
      (\a+0.15,-0.5) -- (\b-0.15,-0.5)
      node[midway,below=5pt] {\small interface \k};
  \node[anchor=east] at (-0.8,2.6) {\small first layer};
  \node[anchor=east] at (-0.8,-1.6) {\small second layer};
\end{tikzpicture}
\caption{The two layers of Subsection~\ref{phys:splitsection}, drawn
for $b=8$ and $m=35$. Above the sites: the chunks of the first layer.
The last chunk has $11$ sites because it absorbs the residue of three
sites. Below: the intervals of the second layer. Each joins the right
half of a chunk to the left half of the next chunk. The alternating
projections on chunk~$i$ and on interface~$i$ are $P_i$ and $Q_i$.}
\label{fig:layers}
\end{figure}

Consider an interval of length $m>2b$, where $b\ge8$, under the uniform law
on one alternating coset of its labels: the set of labels on the
interval is fixed, and their arrangement is uniform among the
arrangements of one parity. In this subsection $\E$, $\Ent$ and
$\|\cdot\|_2$ refer to this law. Form consecutive chunks of
size $b$, merging a final residue into its predecessor. These chunks,
indexed consecutively, are the first layer. At each interface join
the right half of one chunk, rounded down, to the left half of the
next, rounded up, to form the second layer. The intervals within a
layer are disjoint, their lengths lie in $[b/2,2b]$, and each
meets at most two intervals of the other layer; see
Figure~\ref{fig:layers}. The interval
of the second layer that joins chunks $i$ and $i+1$ is called
interface~$i$. Let $P_i$ and $Q_i$ be the alternating projections on
chunk~$i$ and on interface~$i$, that is, the averages $B_I^+$ of
Subsection~\ref{phys:gapsection} over the alternating groups of these
intervals. Their products within the two layers,
\[
 B_1=\prod_iP_i,\qquad B_2=\prod_iQ_i,
\]
are the operators $B_{\mathcal I}^+$ of that subsection for the first
and for the second layer. We write
\begin{equation}\label{phys:def-H}
\begin{aligned}
 S(g)&=\sum_i\|(I-P_i)g\|_2^2+\sum_i\|(I-Q_i)g\|_2^2,\\
 H(g)&=\|(I-B_1)g\|_2^2+\|(I-B_2)g\|_2^2.
\end{aligned}
\end{equation}
Finally, $\Ent_X$ denotes the entropy of a function of a random
variable $X$ under the law of $X$.

\begin{lemma}\label{phys:split}
If $X$ is the set of labels in the first $\lfloor m/2\rfloor$ sites,
then for $f\ge0$,
\[
 \Ent_X\bigl(\E[f^2\mid X]\bigr)\le
 C\bigl((m/b)^2+m\bigr)H(f).
\]
\end{lemma}
\begin{proof}
\step{Step 1: Smoothing by the two layers.}
For~$v=f-\E f$ let~$g=B_2B_1v$ and~$L=\|v\|_2^2-\|g\|_2^2$. The
Pythagorean identity for each projection and the triangle inequality
give~$L+\|v-g\|_2^2\le C H(f)$. Apply the projections one at a time,
first the~$P_i$ and then the~$Q_i$. The iterates~$g_0,g_1,g_2,\dots$
start at~$g_0=v$ and end at~$g$. Fix one of
these projections, and call the component of a function orthogonal to
its range the residual of the function. Immediately after the
projection is applied, the residual of the current iterate is zero. A
later projection that commutes with it does not increase the norm of
the residual. A later projection on an overlapping interval, applied
at step~$j$, increases this norm by at most~$\|g_j-g_{j-1}\|_2$.
There are at most two overlaps in either direction. Since~$S(g)$ is
the sum of the squared norms of the residuals of~$g$,
\begin{equation*}
 S(g)\le4\sum_j\|g_j-g_{j-1}\|_2^2=4L\le C H(f)\,.
\end{equation*}
Thus the passage from~$v$ to~$g$ bounds the sum of the individual
projection losses by the simultaneous layer losses, with squared
displacement at most~$C H(f)$.

The projections preserve the mean, so~$g$ is centered. Apply
Proposition~\ref{smooth:estimate} to~$g$, which plays the role
of~$v$ there, with output~$u=B_1B_2g$. The proposition
bounds~$\Ent(u^2)+\|u\|_2^2$ by~$C(m/b)^2S(g)$ and the second squared
displacement~$\|g-u\|_2^2$ by~$CS(g)$. Put~$w=v-u$. Combining the two
displacement estimates
through~$\|w\|_2^2\le2\|v-g\|_2^2+2\|g-u\|_2^2$ gives
\begin{equation}\label{phys:smoothed}
 \Ent(u^2)+\|u\|_2^2\le C(m/b)^2 H(f),\qquad
 \|w\|_2^2\le C H(f).
\end{equation}

\step{Step 2: A centering inequality.}
We use the following centering inequality for nonnegative $\psi$:
\begin{equation}\label{phys:center}
 \Ent(\psi^2)\le C\left(\Ent((\psi-\E \psi)^2)+\Var(\psi)\right).
\end{equation}
To see this, normalize $\E \psi^2=1$, and put $a=\E \psi$ and
$\sigma^2=\Var(\psi)>0$.
The function $t^2\log t^2-t^2+1$ is at most $4(t-1)^2$ on $[0,2]$
by a Taylor upper bound. For $t>2$, writing $z=t-a$ gives $z>1$ and
$t\le2z$, and the same function is at most
$C z^2(1+\log_+(z^2/\sigma^2))$. Here $\log_+t=\max(0,\log t)$.
Average at $t=\psi$, use $(1-a)^2\le\sigma^2$, and for $z=\psi-a$ observe
$\E z^2\log_+(z^2/\sigma^2)\le\Ent(z^2)+\sigma^2/e$. This proves
\eqref{phys:center},
including the constant case by direct evaluation.

\step{Step 3: Label-set entropy.}
Put $\psi(X)=(\E[f^2\mid X])^{1/2}$. Since $\E f\ge0$, the conditional
$L^2$ triangle inequality gives $|\psi-\E f|^2\le\E[v^2\mid X]$, so
$\Var(\psi)\le\|v\|_2^2$ and
$(\psi-\E \psi)^2\le2\E[v^2\mid X]+2\|v\|_2^2$.
By the variational formula, the entropy of a nonnegative function of
$X$ is the supremum of its integral against the functions
$\phi=\phi(X)$ with $\E e^\phi\le1$; we call such a $\phi$ an entropy
test.
The set $X$ is uniform on a state space with
$\binom m{\lfloor m/2\rfloor}$ elements. Hence the positive part of
every entropy test $\phi$
satisfies $\phi_+\le\log\binom m{\lfloor m/2\rfloor}\le m\log2$
and $\E e^{\phi_+}\le2$.

Fix an entropy test $\phi$. Since $\phi\le\phi_+$ and $\phi_+$ is a
function of $X$, the bound on $(\psi-\E\psi)^2$ and $v=u+w$ give
\[
 \E[(\psi-\E \psi)^2\phi]\le2\E[v^2\phi_+]+2\|v\|_2^2\,\E\phi_+
 \le4\E[u^2\phi_+]+4\E[w^2\phi_+]+2(\log2)\|v\|_2^2 .
\]
Here $\E\phi_+\le\log\E e^{\phi_+}\le\log2$ by Jensen's inequality.
For the remainder $w$ we use the pointwise bound, which gives
$\E[w^2\phi_+]\le m(\log2)\|w\|_2^2$. For $u$ we use the entropy
variational formula under the law of the arrangement. The function
$\phi_+-\log2$ has exponential mean at most one, so
$\E[u^2\phi_+]\le\Ent(u^2)+(\log2)\|u\|_2^2$. We obtain
\[
 \E[(\psi-\E \psi)^2\phi]\le
 4\Ent(u^2)+4(\log2)\|u\|_2^2
 +4m(\log2)\|w\|_2^2+2(\log2)\|v\|_2^2 .
\]
Take the supremum over $\phi$ and apply \eqref{phys:center} and
\eqref{phys:smoothed}, noting that
$\|v\|_2^2\le2\|u\|_2^2+2\|w\|_2^2$. Since $\psi^2=\E[f^2\mid X]$,
this is the assertion.\qedhere
\end{proof}

\subsection{Proof of the comparison theorem}\label{phys:proofsection}

\begin{proof}[Proof of Theorem~\ref{thm:physical}]
\step{Step 1: Entropy decomposition.}
Take $f=\sqrt h$ and fix a row $y$. The entropy loss of $h$ under the
line average of this row is the mean, over the board outside the row,
of the entropy of $h$ under the conditional law of the arrangement of
the row. We first bound this entropy on each alternating coset of
the row, conditionally on the board outside the row.

Starting from the row, split every interval of length greater than
$m_0=128$ into two halves, and repeat with the halves. The intervals
obtained in this way are the
nodes of a binary tree whose root is the row (Figure~\ref{fig:tree}).
An internal node $I$ of size $m$ has two children: $I'$, which
consists of its first $\lfloor m/2\rfloor$ sites, and $I''$, which
consists of the others. The leaves are the nodes of size at most
$m_0$.

\begin{figure}[ht]
\centering
\begin{tikzpicture}[x=0.72cm,y=0.72cm]
  \fill[black!14] (0,0) rectangle (16,0.5);
  \fill[black!14] (8.1,-1.2) rectangle (16,-0.7);
  \fill[black!14] (8.1,-2.4) rectangle (11.95,-1.9);
  \fill[black!14] (10.05,-4.4) rectangle (10.95,-3.9);
  \draw (0,0) rectangle (16,0.5);
  \foreach \a/\b in {0/7.9,8.1/16} \draw (\a,-1.2) rectangle (\b,-0.7);
  \foreach \a/\b in {0/3.9,4.05/7.9,8.1/11.95,12.1/16}
    \draw (\a,-2.4) rectangle (\b,-1.9);
  \node at (4,-3.05) {$\vdots$};
  \node at (13,-3.05) {$\vdots$};
  \foreach \k in {0,...,15} \draw (\k+0.05,-4.4) rectangle (\k+0.95,-3.9);
  \draw[dashed] (10.5,0.9) -- (10.5,-4.8);
  \node[above] at (10.5,0.9) {\small column $j$};
  \node[anchor=east] at (-0.4,0.25) {\small level $0$, the row};
  \node[anchor=east] at (-0.4,-0.95) {\small level $1$};
  \node[anchor=east] at (-0.4,-2.15) {\small level $2$};
  \node[anchor=east] at (-0.4,-4.15) {\small leaves};
\end{tikzpicture}
\caption{The recursion in the proof of Theorem~\ref{thm:physical}.
Each bar is an interval of the row, and an interval longer than $m_0$
is split into two halves. The intervals of level $a$ have length at
most $\lceil n/2^a\rceil$. An internal node contributes a label-set
term and a parity term, and a leaf contributes the entropy of its
internal order. A column $j$ (dashed) meets one interval of each level
(shaded), so its energy $e_{C_j}(f)$ enters only along one branch.}
\label{fig:tree}
\end{figure}
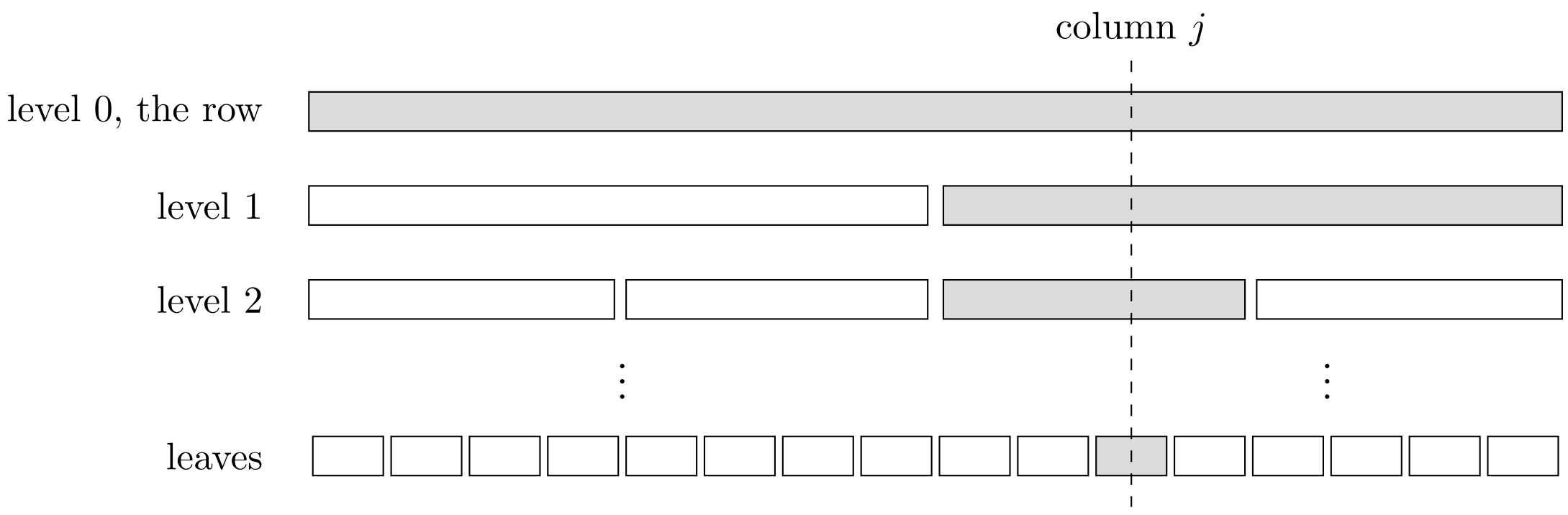

For a node $I$, condition on the board outside $I$. The labels on $I$
are then fixed, and their arrangement is uniform on one alternating
coset. Write $\E_I$ and $\Ent_I$ for the expectation and the entropy
under this conditional law. Both are functions of the board outside
$I$.

Let $I$ be an internal node, let $X$ be the set of labels on $I'$, and
let $\varepsilon$ be the parity of their arrangement on $I'$.
Conditional on $X$,
the remaining law is the equal mixture of two products of alternating
child cosets, indexed by $\varepsilon$. The entropy
chain rule, applied first to $X$ and then to $\varepsilon$, and the
tensorization of entropy for the product law give
\begin{equation*}
 \Ent_I(f^2)\le\Ent_X\bigl(\E_I[f^2\mid X]\bigr)
 +\E_I\bigl[\Ent_\varepsilon\bigl(\E_I[f^2\mid X,\varepsilon]\bigr)\bigr]
 +\E_I\bigl[\Ent_{I'}(f^2)+\Ent_{I''}(f^2)\bigr].
\end{equation*}
Here $\Ent_\varepsilon$ is the entropy under the law of $\varepsilon$
given $X$, which is uniform on the two parities. The first term is the
label-set entropy of Lemma~\ref{phys:split}. The second is the entropy lost
on averaging over the parity bit. The last two are the conditional
entropies on the children, each with coefficient one.

Iterating this decomposition from the root gives, on each alternating
coset of the row,
\begin{equation*}
 \Ent_{\rm row}(f^2)\le\sum_{I\ {\rm internal}}
 \E_{\rm row}\Bigl[\Ent_X\bigl(\E_I[f^2\mid X]\bigr)
 +\Ent_\varepsilon\bigl(\E_I[f^2\mid X,\varepsilon]\bigr)\Bigr]
 +\sum_{I\ {\rm leaf}}\E_{\rm row}\bigl[\Ent_I(f^2)\bigr].
\end{equation*}
Here $\E_{\rm row}$ and $\Ent_{\rm row}$ are $\E_I$ and $\Ent_I$ for
the root, and $X$ and $\varepsilon$ are those of the node $I$. The
iteration leaves a label-set term and a parity term
at each internal node, together with the internal-order entropies at
the leaves.

If $n\le m_0$ no division takes place. The row is then a single leaf,
and the proof uses only the leaf estimate of Step~3, together with the
parity argument at the end of Step~4 when $n$ is even. For $n=2$ the
alternating group of the row is trivial and the parity argument alone
gives the bound.

\step{Step 2: Label-set terms.}
At an internal node $I$ of size $m$ take $b=\lfloor\sqrt m\rfloor$ in
Lemma~\ref{phys:split}. Applied under the conditional law of the node,
the lemma bounds the label-set entropy $\Ent_X(\E_I[f^2\mid X])$ by
$C((m/b)^2+m)H(f)$, where $H(f)$ is the quantity of
\eqref{phys:def-H} for this conditional law.
For each fixed node, integrate the losses over the board outside its
interval and over all parity and label data used to reach it. Its interval
permutations preserve that conditional law, so the resulting projection
losses are exactly their norms under $\pi$:
\begin{equation*}
 \E_\pi[H(f)]=\|(I-B_1)f\|_{2,\pi}^2+\|(I-B_2)f\|_{2,\pi}^2 .
\end{equation*}
On the right-hand side, $B_1$ and $B_2$ are the averages over the two
layers of
short intervals of the node, acting on functions of the whole board.
We can now apply Lemma~\ref{phys:routing} under $\pi$, to
each of the two layers.
Its intermediate physical rotations may leave the conditional space
at a node; this is why we integrate the conditional losses first.
Since $(m/b)^2+m\le Cm$ and $b^2\le m$,
the integrated label-set entropy at $I$ is at most
\begin{equation}\label{phys:node}
 C\left(m^2e_{R_y}(f)+m\sum_{j\in I}e_{C_j}(f)\right).
\end{equation}

\step{Step 3: Parity and terminal intervals.}
To couple the two parity values at a division, transpose the two sites
nearest the cut on each side simultaneously. This involution $s$ preserves
$X$ and sends each conditional coset product uniformly onto the other.
Let $\Ent(a^2,b^2)$ denote the entropy of the function with the two
values $a^2$ and $b^2$ under uniform weights. Apply the two-point
bound $\Ent(a^2,b^2)\le C(a-b)^2$ to the conditional $L^2$ norms of
$f$ on those products. Their
squared difference is at most the conditional squared increment under
$s$. After averaging, the contribution is thus at most $C e_s(f)$.
Since $s$ belongs to the alternating group on the four-site interval
at the cut, this is at most four times $C$ times the corresponding
alternating projection loss. A leaf contributes at most
$C_{m_0}$ times its alternating projection loss by \eqref{phys:center}
and $\Ent(g^2)\le\log(m_0!)\|g\|_2^2$ for its uniform law.

The parity and leaf terms cost at most
$C_{m_0}(e_{R_y}(f)+\sum_{j\in I}e_{C_j}(f))$. Indeed, on an interval
of length between three and $m_0$, every even permutation is a product
of a bounded number of three-cycles. The construction in the proof of
Lemma~\ref{phys:routing} implements each such three-cycle with strides
at most $m_0$ and columns in that interval. Intervals of length at
most two have trivial alternating group.

\step{Step 4: Summation over scales.}
We sum \eqref{phys:node} over all nodes for the row term and over
the branch containing a fixed column for its column term.
At level $a$ of the binary tree, the intervals have total size at
most $n$ and maximum
size $\lceil n/2^a\rceil$. Summing over levels therefore gives
$\sum_I |I|^2\le Cn^2$. Each fixed
column $j$ lies in a single chain of intervals with $\sum_{I\ni j}|I|
\le Cn$. Adding \eqref{phys:node} and the parity and leaf terms bounds
the alternating entropy loss on the row, averaged over its cosets, by
\[
 C\left(n^2e_{R_y}(f)+n\sum_j e_{C_j}(f)\right).
\]
For even $n$ the line group of the row is symmetric. Its extra parity
bit is
coupled by the odd unit rotation $R_y$ and costs at most $C e_{R_y}(f)$ by
the same two-point inequality. For odd $n$ the alternating estimate
already applies on each communicating class. Averaging the row
estimates and the transposed column estimates over the $2n$ lines and
applying \eqref{phys:normalization} proves the theorem.
\end{proof}

\begin{appendices}
\section{Smoothing on a line by two layers of~local~permutations}\label{smooth:section}

We prove the smoothing estimate used in Lemma~\ref{phys:split}.
Averaging over short intervals produces a function whose entropy is
controlled by the local projection losses, with a small $L^2$ error.
The intervals form two layers: consecutive chunks and intervals
straddling their interfaces.

Fix the labels on an interval of length $m$ and restrict their
arrangements to one alternating coset, that is, to one total
permutation parity relative to a reference arrangement. Until we
explicitly pass to the symmetric group, all expectations in this
appendix use the uniform law on this coset. We use the chunks and
interfaces defined before Lemma~\ref{phys:split}, with parameter
$b\ge8$ and $m>2b$, and the notation introduced there. Thus $P_i$ is
the alternating projection on chunk $i$, and $Q_i$ is the alternating
projection on its interface with chunk $i+1$. The operators $B_1$ and
$B_2$ are the products of these commuting projections within each
layer, and $S(v)$ is the sum of the individual losses:
\[
 B_1=\prod_iP_i,\qquad B_2=\prod_iQ_i,\qquad
 S(v)=\sum_i\|(I-P_i)v\|_2^2+\sum_i\|(I-Q_i)v\|_2^2.
\]

\begin{proposition}[Smoothing by two layers]\label{smooth:estimate}
For every centered $v$, the function $u=B_1B_2v$ satisfies
\[
 \|v-u\|_2^2\le C S(v),\qquad
 \Ent(u^2)+\|u\|_2^2\le C(m/b)^2 S(v).
\]
\end{proposition}

The order $B_1B_2v$ applies the interface averages before the chunk
averages. Chunk averaging removes internal order within each chunk
parity class; paired odd permutations in neighboring chunks remove
the remaining dependence on the separate chunk parities. The proof
has three parts.
Lemma~\ref{smooth:local} in Subsection~\ref{smooth:two-chunk}
establishes this local reduction to
label sets and the displacement bound.
Lemma~\ref{smooth:cost} in Subsection~\ref{smooth:exchange-cost}
bounds the cost of exchanging
labels between neighboring chunks by $C/b$ times the corresponding
projection loss. Subsection~\ref{smooth:chunk-entropy} completes the
proof of Proposition~\ref{smooth:estimate}: we extend the resulting
function of chunk label sets to the symmetric group and use a
binary-tree entropy decomposition to obtain the stated factor
$(m/b)^2$. The letters $p$, $q$ and $\lambda_j$ have meanings local to
this appendix. They are unrelated to the exponent $p$ and to the
shapes $\lambda$ of Section~\ref{glob:section}.

\subsection{The two-chunk transition}\label{smooth:two-chunk}

In this subsection and the next, $v$ is an arbitrary function on the
coset, not necessarily centered, and $u=B_1B_2v$.
Fix an interface $i$ and put
\begin{equation*}
 v_i=\left(\prod_{j\notin\{i,i+1\}}P_j\right)\left(\prod_{h\ne i}Q_h\right)v,\qquad T_i=P_iP_{i+1}Q_i\,.
\end{equation*}
Let $l$ and $r$ be the sizes of chunks $i$ and $i+1$, and
$p$ and $q$ the numbers of left and right sites used by the interface.
Thus $p/l,q/r\in[1/3,2/3]$ and each of $p$ and $q$ is at least two.
We call the arrangement outside these two chunks the exterior.
When we condition on the exterior,
the set of labels on the two chunks and the total parity of their
arrangement are fixed, so this arrangement is uniform on an
alternating coset. We then write $A$ for the set of
the $l$ labels in chunk~$i$, the chunk on the left of the interface.

\begin{lemma}[Two-chunk transition]\label{smooth:local}
For every interface $i$ the following hold.
\begin{enumerate}
\item[(i)] $u=T_iv_i$, and the local losses
$d_i=\|v_i\|_2^2-\|u\|_2^2$ satisfy
\begin{equation}\label{smooth:loss}
 d_i\le \|(I-Q_i)v\|_2^2+
             \sum_{j=i,i+1}\|(I-P_j)B_2v\|_2^2,\qquad
 \sum_i d_i\le C S(v).
\end{equation}
Also $\|v-u\|_2\le\sqrt{2S(v)}$.
\item[(ii)] Every function in the range of $T_i$ is invariant under
the even permutations of the sites of chunks $i$ and $i+1$ that map
each of the two chunks to itself. Conditional on the arrangement
outside these two chunks, it therefore depends only on the set $A$.
\item[(iii)] Under the same conditioning, $T_iT_i^*$ acts on functions
of $A$ as the transition operator of the following random step. Choose
independently a uniform $p$-subset of $A$ and a uniform $q$-subset of
the labels in chunk $i+1$. Redistribute their union uniformly, placing
$p$ of its labels in chunk~$i$.
\end{enumerate}
\end{lemma}

\begin{proof}
\step{Part (i).}
The members of either layer commute, and the projections on the chunks other than $i$ and $i+1$ commute with $Q_i$. Consequently $u=T_i v_i$.
Successive orthogonal projection gives \eqref{smooth:loss}.
For the first inequality, use $Q_iv_i=\left(\prod_{j\notin\{i,i+1\}}P_j\right)B_2v$ and discard commuting
contractions in each squared increment. For the second, when bounding
$\|(I-P_j)B_2v\|_2$, first drop the projections $Q_h$ whose interfaces
are disjoint from chunk $j$.
At most two remain. Telescoping their product shows that the
norm is at most $\|(I-P_j)v\|_2+\sum_{h\in\{j-1,j\}}\|(I-Q_h)v\|_2$,
with missing boundary interfaces omitted. Each summand occurs only a bounded number of times.
Commutation within each layer also gives the displacement bound
$\|v-B_1B_2v\|_2\le\|(I-B_1)v\|_2+\|(I-B_2)v\|_2\le\sqrt{2S(v)}$.

\step{Part (ii).}
We show that the output of $T_i$ depends only on the labels on
the left of the interface.
Condition on the exact arrangement outside the two chunks.
For a fixed
set of labels on the left, the projection $P_iP_{i+1}$ removes
internal order within each pair of chunk parity classes. The fixed
total parity links those classes, leaving one possible parity bit.

Choose in each chunk a transposition of two of its sites used by the
interface, which is possible because $p,q\ge2$, and let
$s$ be their product. This permutation reverses both chunk parities
and preserves their total parity. It belongs to the interface
alternating group and normalizes the subgroup underlying
$P_iP_{i+1}$, so $U_s$ commutes with $P_iP_{i+1}$ and satisfies
$U_sQ_i=Q_i$. Consequently $U_sT_i=T_i$: the output cannot
distinguish the two values of the remaining parity bit. Together with
the two chunk alternating groups, $s$ generates all within-chunk
permutations having even total parity, since an odd permutation in
each chunk becomes even in each after multiplication by $s$. The
output of $T_i$ therefore depends only on the $l$-set $A$ of labels
on the left, conditional on the exterior; it retains neither the
internal orders nor the separate chunk parities.

\step{Part (iii).}
We have $T_iT_i^*=P_iP_{i+1}Q_iP_iP_{i+1}$. On functions of $A$ this
operator acts as the stated transition. Indeed, each
chunk alternating group acts transitively on subsets of the prescribed
size, as does the alternating group of the interface union. Any
necessary parity correction for a subset image can be made by an
internal transposition preserving that subset. Thus the stated
transition law holds within the given alternating coset.
\end{proof}

\subsection{Spectral estimate and exchange cost}\label{smooth:exchange-cost}

We compare the subset transition of Lemma~\ref{smooth:local}(iii)
with the complete transposition
Laplacian. Both operators are diagonal on the same decomposition of
subset functions by degree. Fix an interface $i$ and condition on the
arrangement outside chunks $i$ and $i+1$, as in
Lemma~\ref{smooth:local}. Interchanging the chunks if necessary, we
may suppose that $l\le r$. On $l$-sets let $V_j$ be the
span of the monomials $\one_{\{D\subset A\}}$ with
$|D|=j$. These spaces are nested, as summing a monomial over all
one-point extensions of $D$ multiplies it by $l-j$. The space $V_l$ is the
entire function space. Set $W_j=V_j\ominus V_{j-1}$, with $V_{-1}=0$.

The second object of this subsection is the cost of exchanging one
label between the two chunks. By Lemma~\ref{smooth:local}(ii) the
output $u$ is a function of $A$. Fix a site of chunk $i$ and a site of
chunk $i+1$, and exchange the labels at these two sites.
At this stage the
exchange acts on the left label set: it replaces one left label by
one right label. This is well defined for the subset function $u$,
although the corresponding site transposition leaves
the original alternating coset.
Subsection~\ref{smooth:chunk-entropy} identifies this subset
increment with an ordinary site-transposition increment of the
extended function.
Let $\delta_i$ be the mean of the squared increment of $u$ under this
exchange, over the arrangement in the two chunks and over the
exterior.

\begin{lemma}[Spectrum and exchange cost]\label{smooth:cost}
For every interface $i$ the following hold.
\begin{enumerate}
\item[(i)] The operator $T_iT_i^*$ of Lemma~\ref{smooth:local}(iii)
acts on each $W_j$ as a scalar $\lambda_j$, and
\begin{equation}\label{smooth:eigen}
 0\le\lambda_j\le (1-pq/(lr))^j\le(8/9)^j,\qquad 1\le j\le l.
\end{equation}
\item[(ii)] The cost of one exchange satisfies
\begin{equation}\label{smooth:interface}
 \delta_i\le C d_i/b .
\end{equation}
\end{enumerate}
\end{lemma}

\begin{proof}
\step{Part (i): Scalar action.}
An operator on subsets commuting with every relabeling sends a fixed
degree-$j$ monomial to a function of $|A\cap D|$, the same function
for each $D$ of size $j$. Expanding in
$\binom{|A\cap D|}{a}$ with $0\le a\le j$ shows that it acts by a scalar on
$V_j/V_{j-1}$. For a self-adjoint operator it acts by that scalar on
$W_j$. Let $\lambda_j$ be this scalar for $T_iT_i^*$.

\step{Part (i): Eigenvalues of the subset transition.}
To bound $\lambda_j$, take disjoint pairs of labels $(a_h,a'_h)$ with
$1\le h\le j$, and define
\[
 \Phi_j(A)=\prod_{h=1}^j
       \left(\one_{\{a_h\in A\}}-\one_{\{a'_h\in A\}}\right).
\]
Each monomial of degree below $j$ misses both labels of some pair, so
swapping that pair proves $\Phi_j\in W_j$. At a set containing all
$a_h$ and no $a'_h$ we have $\Phi_j=1$. In the subset transition, if
the sampled union contains both partners of some pair, its conditional
mean after redistribution is zero by exchangeability. Otherwise the
absolute value of $\Phi_j$ is at most one. If $\nu$ of the labels $a_h$
were sampled on the left, the probability that the right sample avoids
all their
partners is at most $(1-q/r)^\nu$. For $0\le z\le1$,
\[
 \E z^\nu\le \bigl(1-(1-z)p/l\bigr)^j.
\]
Indeed, in the expansion of $\prod_h[z+(1-z)\one_{\{a_h\ {\rm unsampled}\}}]$,
each event of avoiding $w$ prescribed left labels has probability at most
$(1-p/l)^w$.
Evaluation at the chosen set and positivity of $T_iT_i^*$ yield
\eqref{smooth:eigen}.

\step{Part (ii): Energy of the output.}
The complete transposition Laplacian on $l$-sets, with each unordered
pair of labels swapped at rate one, acts on $W_j$ by
$j(l+r-j+1)$. To verify the scalar, evaluate it on $\Phi_j$ at the same
initial set. Swapping partners contributes $2j$. Swapping opposite-side
labels from distinct marked pairs contributes $j(j-1)$, while a marked
and an unmarked label on opposite sides contribute
$j(r-j)+j(l-j)$. All other swaps preserve the value.

The singular values of $T_i$ relate this energy to the loss $d_i$.
Fix the exterior and write $u=\sum_ju_j$ with $u_j\in W_j$.
In a singular-value decomposition, an output component of squared
norm $\|u_j\|_2^2$ in $W_j$ with $\lambda_j>0$ requires an input
of squared norm at least $\lambda_j^{-1}\|u_j\|_2^2$. Every output in a
zero eigenspace vanishes. Thus the local loss $\|v_i\|_2^2-\|u\|_2^2$, with the
current exterior fixed, bounds $\sum_{j:\lambda_j>0}
(\lambda_j^{-1}-1)\|u_j\|_2^2$. The energy of $u$ is
$\sum_jj(l+r-j+1)\|u_j\|_2^2$, which is at most
$(l+r)\sum_jj\|u_j\|_2^2$. By \eqref{smooth:eigen},
$j\lambda_j\le C(1-\lambda_j)$, that is, $j\le C(\lambda_j^{-1}-1)$
when $\lambda_j>0$, and $l+r\le4b$. The complete
transposition energy of the output is therefore at most $Cb$ times
the local loss.

\step{Part (ii): Cost of one exchange.}
At any $A$
only the $lr$ label pairs in $A\times A^c$ contribute to this energy,
where $A^c$ is the set of labels in chunk $i+1$.
Conditional on $A$, the two labels at any prescribed
cross-chunk pair of sites are uniform on $A\times A^c$ by the separate
transitivity of the two chunk alternating groups. In particular
$\delta_i$ does not depend on the choice of the two sites. The energy is
therefore $lr/2$ times the squared increment of a single fixed
cross-chunk exchange, averaged over the two chunks.
As $lr\ge b^2$, integration over the exterior gives
\eqref{smooth:interface}.
\end{proof}

\subsection{Entropy of the chunk label sets}\label{smooth:chunk-entropy}

\begin{proof}[Proof of Proposition~\ref{smooth:estimate}]
\step{Step 1: The output depends only on the chunk label sets.}
By Lemma~\ref{smooth:local}, $u=T_i v_i$ at every interface $i$, so
$u$ is invariant under an
odd permutation in each of two neighboring chunks, as well as under
every chunk alternating group. These paired odd permutations generate
every even change of chunk parities: applying neighboring pairs along
a path reverses the endpoint parities and changes each intermediate
parity twice. Two arrangements in the original coset with the same
ordered chunk label sets differ by within-chunk permutations of even
total parity. The paired flips match their chunk parities, and the
chunk alternating groups then match their internal orders. Hence $u$
depends only on the ordered list of chunk label sets.

\step{Step 2: Extension to the symmetric group.}
Extend $u$ to the full symmetric group $S_m$ by assigning every
arrangement the value determined by its ordered chunk label sets. A
transposition within any one chunk interchanges the two global parity
classes while preserving all these sets. Their distribution is
therefore the same in either class and under the uniform symmetric
law, so the extension preserves~$\E u$, $\|u\|_2^2$
and~$\Ent(u^2)$. Conditional on the chunk label sets, the labels at a
prescribed pair of sites in chunks $i$ and $i+1$ are uniform on the
product of the two sets under either law, by separate transitivity of
the chunk alternating groups. Thus the mean squared
site-transposition increment of the extension is precisely
$\delta_i$. We now use the uniform law on $S_m$ for $u$ and its
increments; the quantities $d_i$ and $S(v)$ still refer to the
original alternating coset. Ordinary site transpositions are now
defined on the whole space and preserve its uniform law.

\step{Step 3: A binary tree of chunks.}
To estimate the entropy of the extension, form a binary tree of
lists of chunks. Bisect the list of all chunks, and then each
resulting list of more than one chunk, into two lists whose numbers of
chunks differ by at most one. This tree divides lists of whole
chunks, unlike the divisions of site intervals in the proof of
Theorem~\ref{thm:physical}; its leaves are the original chunks.
Identify each node $I$ with the union of the sites in its chunks, so
$|I|$ counts sites. For each internal node $I$, the two child sizes
as fractions of $|I|$ lie in $[1/5,4/5]$ because each chunk has
length between $b$ and $2b$. Under the symmetric
law, conditional on the arrangement outside $I$, the labels in the
first child form
a uniform subset of fixed cardinality. If $\xi$ is that subset, put
$\psi(\xi)=(\E[u^2\mid\xi])^{1/2}$ with the arrangement outside $I$
fixed.
Swapping two labels couples the conditional laws at the corresponding
subsets. The difference of their conditional $L^2$
norms is at most the $L^2$ norm of the coupled increment of $u$.
The logarithmic Sobolev estimate of Lee and Yau for the Bernoulli--Laplace
model \cite{leeyau},
with each unordered transposition at rate one, bounds the subset
entropy of $\psi^2$ by $C/|I|$ times the sum of squared increments over
label pairs for these balanced fractions. This normalization appears in
\cite[Eq.~(2.66)]{caputosohier}.
For each assignment the sum over label pairs equals the sum over site
pairs. Conditional on the child label sets the two child permutations
are independent, so entropy tensorization and induction give
\begin{equation}\label{smooth:tree}
 \Ent(u^2)\le C\sum_{I\ {\rm internal}}\frac1{|I|}
       \sum_{\{x,y\}\subset I}
       \E\!\left[(u((x\,y)\sigma)-u(\sigma))^2\right].
\end{equation}
Leaves contribute zero since $u$ is constant under permutations of
their labels.

\step{Step 4: Reduction to neighboring chunks.}
We reduce the exchanges in \eqref{smooth:tree} to neighboring chunks.
For a pair of sites in different chunks of $I$, choose a path from
one to the other with one site in each intervening chunk. The endpoint
transposition is a product of adjacent transpositions along the path,
first forward and then back along all but the last edge. The number
of factors is at most $2|I|/b$ and each interface is crossed at most
twice. The uniform law on the symmetric group is preserved by all
prefixes of this word. The Cauchy--Schwarz inequality along the word
and \eqref{smooth:interface} bound the contribution from $I$ on the
right-hand side of
\eqref{smooth:tree} by
\[
 \frac{C|I|^2}{b^2}\sum_{\text{interfaces }i\text{ inside }I}d_i .
\]
Indeed, there are at most $|I|^2$ site pairs and same-chunk pairs have
zero increment. For a fixed interface the sizes of its containing
nodes decrease by a factor at most $4/5$ at each descent. Summing and
using \eqref{smooth:loss} yields $\Ent(u^2)\le C(m/b)^2 S(v)$.

\step{Step 5: The $L^2$ bound.}
Finally, the entropy estimate also gives the $L^2$ bound. The
estimates for projection losses, interface increments and
\eqref{smooth:tree} apply to inputs of arbitrary mean. We may therefore
apply the entropy estimate to $1+t v$, whose output is $1+t u$.
The coefficient of $t^2$ in $\Ent((1+t u)^2)$ is $2\Var(u)$. Since
$\E u=\E v=0$, this bounds $\|u\|_2^2$ as asserted. The displacement
bound is part of Lemma~\ref{smooth:local}(i).
\end{proof}

\section{Estimates for the exceptional configurations}\label{exc:section}

We prove the estimates \eqref{glob:heavy} and \eqref{glob:both} of
Step~6 in the proof of Lemma~\ref{glob:low}. The notation and the
assumptions are those of that proof, and as there a pair means a
designated pair. We use two facts from its Step~4. Given a row orbit
$\mathcal R$, different rows are independent, and in a row with at
least two unmarked sites the marked labels occupy uniformly chosen
distinct columns. Moreover, every event has conditional probability
at most $2^{t/(n-1)}$ times its probability under the symmetric
conditional law, under which the marked labels of every row occupy
uniformly chosen distinct columns. This is the density bound. In each
of the two estimates we prove the first bound. The second bound
follows by interchanging rows and columns.

\begin{proof}[Proof of \eqref{glob:heavy}]
A heavy column contains more than $L$ of the $t$ marked labels, so at
most $u=\lfloor t/L\rfloor$ columns are heavy. By \eqref{glob:counts}
and the assumption $k\le n^{21/20}$, we have $u\le2k/L$ and
$u/n\le n^{-19/100}$ for all sufficiently large $n$. If $u=0$ the
event $\mathcal H_{\rm col}$ is empty. Otherwise, on
$\mathcal H_{\rm col}$ one can find a set of $u$ columns and a set of
$\lceil k/16\rceil$ marked labels such that all these labels lie in
these columns. Fix a set of $u$ columns and a set of
$\lceil k/16\rceil$ marked labels. Under the symmetric conditional
law, the probability that $m$ given labels of one row all lie in the
$u$ given columns is
\[
 \frac un\cdot\frac{u-1}{n-1}\cdots\frac{u-m+1}{n-m+1}
 \le\Bigl(\frac un\Bigr)^{m}.
\]
Different rows are independent, so the probability that all
$\lceil k/16\rceil$ given labels lie in the given columns is at most
$(u/n)^{\lceil k/16\rceil}$. There are at most $n^u$ sets of $u$
columns and at most $2^t$ sets of marked labels. The union bound over
these sets and the density bound give
\[
 \mathbb P(\mathcal H_{\rm col}\mid\mathcal R)
 \le 2^{t/(n-1)}n^u2^t (u/n)^{\lceil k/16\rceil}.
\]
The logarithm of the right-hand side is at most
\[
 \Bigl(\frac{2}{n-1}+2\Bigr)k\log2+\frac{2k}{L}\log n
 -\frac{19}{1600}\,k\log n
 =\Bigl(-\frac{19}{1600}+o(1)\Bigr)k\log n,
\]
where the term $o(1)$ does not depend on $k$. Since $19/1600>1/200$,
this proves the first bound of \eqref{glob:heavy}.
\end{proof}

\begin{proof}[Proof of \eqref{glob:both}]
Fix a row orbit $\mathcal R$. It determines the row of every marked
label, and hence the set of heavy rows.

\step{Step 1: Pairs on light lines.}
On $\mathcal N$ the heavy rows contain at most $k/16$ marked labels,
and so do the heavy columns. At most $k/8$ pairs therefore have a
label on a heavy line, and by \eqref{glob:counts} at least
$s-k/8\ge3k/8$ pairs have all their lines light. Put
$q=\lceil k/4\rceil$. Then at least $q$ pairs have all their lines
light.

\step{Step 2: Candidates and witnesses.}
On $\mathcal N$, a pair whose lines are all light is not free, so its
rectangle meets the rectangle of another cluster. Two rectangles meet
exactly when their row sets meet and their column sets meet.
Accordingly, a \emph{candidate} for a pair consists of a cluster
other than the pair that has a label in a row of the pair, of a label
of the pair, and of a label of that cluster. The candidate is a
\emph{witness} if its two labels lie in the same column. The
candidates of a pair are determined by $\mathcal R$, and on
$\mathcal N$ every pair with light lines has a witness. In
Figure~\ref{fig:exchange}, the singleton with the labels $5$ and $3$
is a witness for the pair $\{4,5\}$.

\step{Step 3: The number of lists.}
A \emph{list} consists of $q$ pairs whose rows are light and of one
candidate for each of them. By Steps~1 and~2, on $\mathcal N$ there
is a list all of whose candidates are witnesses. The $q$ pairs can be
chosen in at most $2^s\le2^k$ ways. The rows of a pair of the list
are light, so they contain at most $2L$ marked labels, and at most
$2L$ clusters have a label in these rows. Each of these clusters
gives at most four candidates. Hence there are at most $2^k(8L)^{q}$
lists.

\step{Step 4: The probability of one list.}
Fix a list. We bound the conditional probability, given $\mathcal R$,
that all its candidates are witnesses. Join the two labels of each
candidate of the list by an edge, and consider the graph formed by
these edges and their endpoints. If all candidates are witnesses,
then all the labels of a connected component of this graph lie in one
column. The labels that the candidates take from the $q$ pairs of the
list are distinct and lie in light rows. So at least $q$ labels in
light rows are endpoints of edges.

A light row has at least two unmarked sites. By the first fact
recalled above, its marked labels occupy uniformly chosen distinct
columns, independently of the other rows, for symmetric and for
alternating line groups. The density bound is therefore not needed
here.

Reveal first the columns of all marked labels in heavy rows,
and then treat the components one by one. If a component contains a
label in a heavy row, the column of this label is already revealed,
and every label of the component in a light row has to take this
column. If all $r$ labels of a component are in light rows, reveal
the column of one of them. The other $r-1$ labels have to take this
column, and $r-1\ge r/2$ because $r\ge2$. Altogether at least $q/2$
labels in light rows have to take a prescribed column.

Given
everything revealed before, the column of a label in a light row is
uniformly distributed over the columns of its row that are not
occupied by labels revealed earlier. At most $L-1$ columns are
occupied in this way, so each prescribed column is taken with
conditional probability at most $(n-L)^{-1}$. The probability that
all candidates of the list are witnesses is therefore at most
$(n-L)^{-q/2}$.

\step{Step 5: Conclusion.}
Summing over the lists gives
\[
 \mathbb P(\mathcal N\mid\mathcal R)
 \le 2^k(8L)^{q}(n-L)^{-q/2}\le n^{-k/40}
\]
for all sufficiently large $n$. For the last inequality, recall that
$q\ge k/4$ and $L\le n^{1/4}$. The logarithm of the middle term is
then at most $(-1/16+o(1))\,k\log n$, where the term $o(1)$ does not
depend on~$k$.
\end{proof}

\end{appendices}

\section*{Use of AI tools}
\phantomsection
\addcontentsline{toc}{section}{Use of AI tools}
The author developed the proof of Theorem~\ref{thm:global} with assistance of GPT-6 Astra (OpenAI), and used Claude Fable 5.2
(Anthropic) to revise and polish the writing of the whole paper and to check
the text.

\end{document}